\documentclass[12pt]{amsart}

\usepackage[all]{xy}
\usepackage{a4wide}
\usepackage{geometry}
\usepackage{amsmath}
\usepackage{amssymb}
\usepackage{amsthm}  
\usepackage{amscd}       
\usepackage{subfigure}
\usepackage{epsfig}
\usepackage{color}
\usepackage{graphicx}      % See geometry.pdf to learn the layout options. There are lots.
\usepackage{graphicx}				% Use pdf, png, jpg, or eps§ with pdflatex; use eps in DVI mode
\usepackage{mathrsfs}
\usepackage{amssymb}
\usepackage{color}
\usepackage{changes}
\usepackage[skip=2pt,font=scriptsize]{caption}
\newcommand{\R}{\mathbb{R}}
\newcommand{\Z}{\mathbb{Z}}

\newcommand{\N}{\mathbb{N}}

\newcommand{\cT}{\mathcal{T}}

\newcommand{\cD}{\mathcal{D}}

\theoremstyle{plain}
\newtheorem{theorem}{Theorem} [section]
\newtheorem{lemma}[theorem]{Lemma}
\newtheorem{prop}[theorem]{Proposition}

\providecommand{\keywords}[1]{\textbf{\textit{Index terms---}} #1}

\numberwithin{equation}{section}

\begin{document}

\title{On the Penrose and Taylor-Socolar Hexagonal Tilings }
\date{}
\author{Jeong-Yup Lee $^{1}$
%\footnote{The first author would like to acknowledge
%that this work was supported by a National Research Foundation
%of Korea (NRF) Grant funded by the Korean
%Government (MSIP) (No. 2014004168) and by research fund of Catholic Kwandong %University(CKURF-201604560001). She is also grateful for the support by the Korea %Institute for Advanced Study (KIAS) . 
, Robert V.\ Moody $^{2}$}

%\address{$^{1}$ Department of Mathematics Education, Catholic Kwandong University,
%Gangneung, Gyeonggi-do, 210-701, Korea, jylee@cku.ac.kr \\
%$^{2}$
%Department of Mathematics and Statistics,
%University of Victoria, Victoria, British Columbia V8W 3P4, Canada, rvmoody@mac.com}

%\corres{jylee@cku.ac.kr, 82-33-649-7776(Tel.), 82-33-642-7716(Fax)}

\maketitle

\begin{abstract} We study the intimate relationship between the Penrose and the Taylor-Socolar tilings, 
within both the context of double hexagon tiles and the algebraic context of hierarchical inverse sequences of triangular lattices. This unified approach produces both types of tilings together, clarifies
their relationship, and offers straightforward proofs of their basic properties. \end{abstract}

\vspace{5mm}

\noindent
\keywords{Keywords: Penrose tiling, Taylor-Socolar tiling, double hexagon tiling, nested triangularizations,
inverse sequences of hierarchical lattices 
}

%\MSC{52C23}

%%%%%%%%%%%%%%%%%%%%%%%
\section{Introduction}\label{intro}
%%%%%%%%%%%%%%%%%%%%%%%

 From the very beginning, aperiodic tilings have played a significant role in unravelling the mysteries of aperiodic crystals. Knowing what is mathematically possible has often turned out to be a crucial element in conceiving what might be physically realizable. In this paper we discuss two remarkable aperiodic tilings of the plane that are built out of one of the most basic of all crystallographic structures: the standard periodic hexagonal lattice.

Also from the very beginning, there arose the question of what might be the minimum number of different prototiles necessary for a system of tiles and corresponding matching rules that permit, and only permit, aperiodic tilings. The very first aperiodic tilings involved thousands of prototiles. The famous aperiodic tilings of the plane like the rhombic Penrose and the Ammann-Beenker tilings are each based on just two prototiles, the allowable motions being translations and rotations. This of course immediately raises the question of whether aperiodic tilings based on just one prototile are possible.

Taylor\cite{Taylor} and Taylor and Socolar\cite{ST} introduced a planar aperiodic tiling which can be built from a single hexagonal prototile allowing translations, rotations, and {\bf reflections}. The tiling is based on the familiar hexagonal tiling of the plane, but if one distinguishes the prototile in its direct and reflected forms, then the matching rules allow only aperiodic tilings to appear. Their work revived interest in much earlier work of R. Penrose. In \cite{Penrose}, he had introduced an aperiodic tiling, also based on marked hexagonal tiles, but additionally involving two other types
of thin edge tiles and small corner tiles, which he called  a $(1+\epsilon+\epsilon^{2})$-tiling. However, already in this paper he had introduced arrowed double hexagon tiles as an alternative way to represent the tiling, 
see Fig.~\ref{doubleArrowing}. Later, in his online notes \cite{PenroseTwistor} he expanded upon the double tile theme, and pointed out the essential matching rules that make them work. 

In both the Taylor-Socolar and Penrose cases, the `matching rules' somewhat stretch the original notion of matching rules, so there can remain some controversy about whether these are strictly aperiodic monotiles. However, that is not an issue here. These tilings are interesting, and puzzling too, for although the Taylor-Socolar hexagonal tilings, henceforth called Taylor-Socolar tilings or T-S tilings, and the Penrose hexagonal tilings (Penrose tilings\footnote{Since Penrose tilings based on five-fold symmetry are so much a part of the aperiodic culture, we should emphasize that the Penrose tilings of this paper are based on hexagons and have nothing to do with the rhombic or kite/dart Penrose tilings.}) seem deeply related, that relationship 
is somewhat obscure. The two tilings are not mutually locally derivable (MLD) \cite{BGG, BGBook} in the technical sense, but are in mutually derivable in a rather different sense that we shall explain.

In \cite{TSMS} we put forward a development of the T-S tilings based on the underlying hierarchical system of nested equilateral triangles that are so prominent in both the T-S tilings and the Penrose tilings. The aperiodicity of the tilings comes from this hierarchical structure, and indeed these tilings seem to have been invented with precisely this feature in mind. The structure of nested triangles has an algebraic interpretation as an inverse system of finite groups, arising from the standard triangular lattice and its natural triangular sublattices, and is closely related to the $2$-adic integers. In \cite{TSMS} we made this algebraic interpretation the basis out of which we constructed the T-S tilings. In fact, as long as the nested system of triangles is generic, meaning that it is free of singularities (like points which are simultaneously the vertices of triangles of unbounded size), then there is a unique T-S tiling belonging to it. We shall see that the same type of mathematics applies to the Penrose tilings, and not surprisingly the two inverse systems are deeply connected. 

A more detached look at the double hexagon tilings reveals that they actually incorporate both types of tilings simultaneously. This association is not entirely new \cite{BGG, BGBook}, but in this paper it is the double hexagon tilings that are taken as the fundamental objects, and they serve as the parents of the two individual types (Penrose and T-S) of hexagonal tilings. Thus we may think of the two tilings as siblings of each other. Algebraically a double hexagon tiling corresponds to a matched pair of inverse sequences, and with these we can see how the algebra and geometry fit seamlessly to elucidate each other.  The paper offers a unified treatment of the two tilings along with proofs of the implied hierarchical structuring and the aperiodicity.

%%%%%%%%%%%%%%%%%%%%%%%%%%%%%%%%
\section{Double hexagon tiles and their tilings}\label{2hexTiles}
%%%%%%%%%%%%%%%%%%%%%%%%%%%%%%%%

\begin{figure}[h]
\centering
\includegraphics[width=0.6\textwidth]{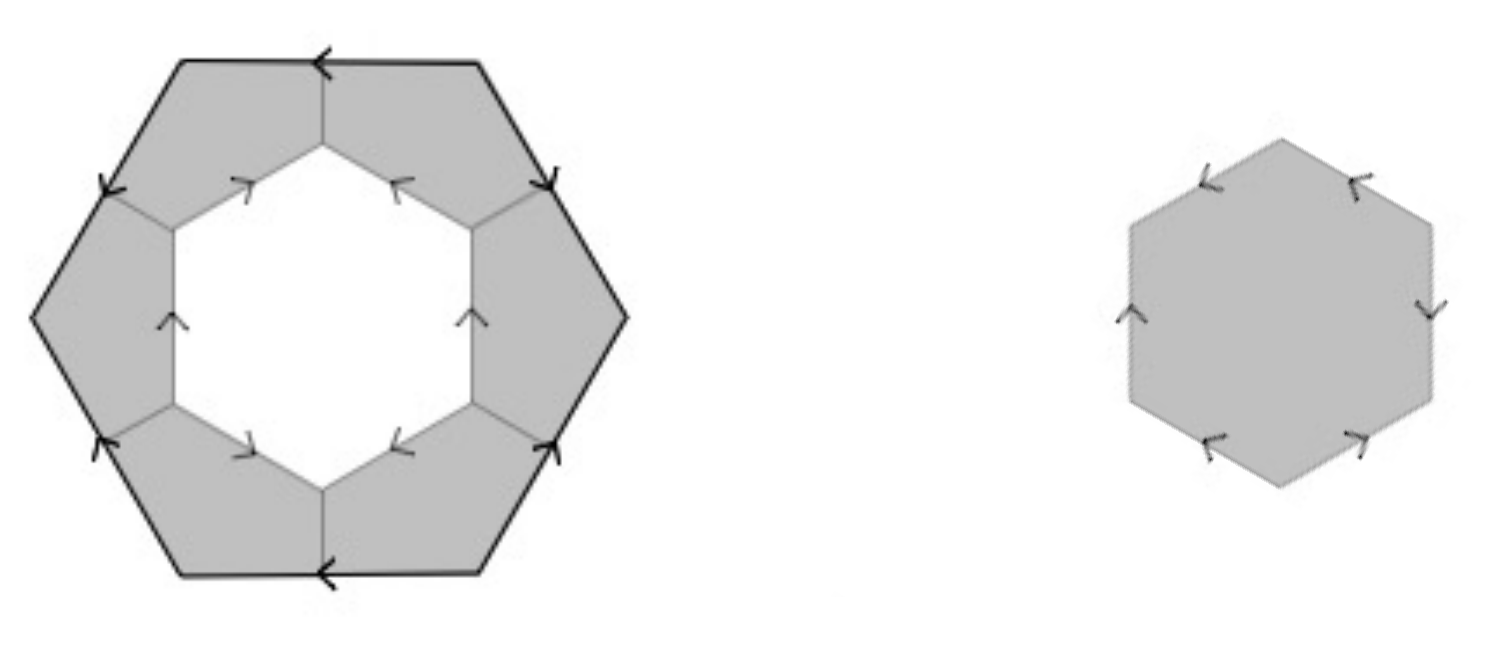}
\caption{A well-arrowed hexagon is shown on the right. It has one pair of opposite sides whose arrows face in the same direction,
thus providing an orientation for the tile.  A double hexagon tile is on the left, the key feature being that the orientations of the inner and outer arrowed hexagons are at right-angles. When three double hexagon tiles meet at a vertex, the gray parts around the vertex form corner hexagons, see Fig.~\ref{doubleHexPatch}. The assumption of {\bf legality} allows completion of the arrowing on the corner hexagons to well-arrowed hexagons, as indicated on the right.}
\label{doubleArrowing}
\end{figure}

 An arrowed hexagon is a regular hexagon in which each side has been given a direction, indicated by an 
arrowhead. An arrowed hexagon is called {\bf well-arrowed} if, up to rotation,  the arrows form the pattern shown
on the right side of Fig.~\ref{doubleArrowing}. In fact all three hexagons in this figure are well-arrowed. The structure of the well-arrowed hexagon gives it a well-defined orientation
in the plane, namely that provided by the two parallel arrows facing in the same direction. 

If we start at any edge of a well-arrowed hexagon and then look every alternate edge as we go around it, we notice that arrows on the three edges always have mixed type of  clockwise and counter-clockwise. 
We notice also the useful
fact that if we have a hexagon and three alternate edges have been arrowed so as to be of mixed type, then there is a unique
way to complete the arrowing to make it into a well-arrowed hexagon. 

A hexagonal tiling of the plane with well-arrowed hexagons is called {\bf well-matched} if the hexagons meet 
edge-to-edge and the arrows of these coinciding edges also coincide -- that is, they point in the same direction. 
We are only interested in well-matched tilings of arrowed hexagons. 

A double hexagon tile (or double hex tile) consists of a pair of well-arrowed hexagons, one within the other, as shown on the left side
of Fig.~\ref{doubleArrowing}. The inner hexagon is centered within the outer one with its orientation at right-angles
to the orientation of the outer one. Its size is chosen so as to make the outer hexagon  $3$ times the area
of the inner hexagon (so there is a linear scaling factor of $\sqrt 3$). There are, up to rotational symmetry, only
two double hexagons (of any particular size), see Fig.~\ref{NewDeco-on-PenroseTile}. When three double hexagons meet at a vertex, the gray parts around the vertex form another hexagon. We call these types of hexagons {\bf corner hexagons}. 

Suppose that we have a hexagonal tiling of the plane with double hexagons. By this we mean
that we are using the outer hexagons as the tiles. 
Suppose that from the perspective of the arrowing on the outer hexagons this hexagonal tiling is well-matched. 
Three outer hexagons meet at every common vertex $v$, and the three edges of the corresponding inner hexagons
that are closest to $v$ form three edges of a corner hexagon $H$ centered on $v$. With the terminology
introduced above we can ask whether or not these three arrows are mixed. If they are mixed then we can extend the arrowing to make $H$ well-arrowed.  

\begin{figure}[h]
\centering
\includegraphics[width=0.7\textwidth]{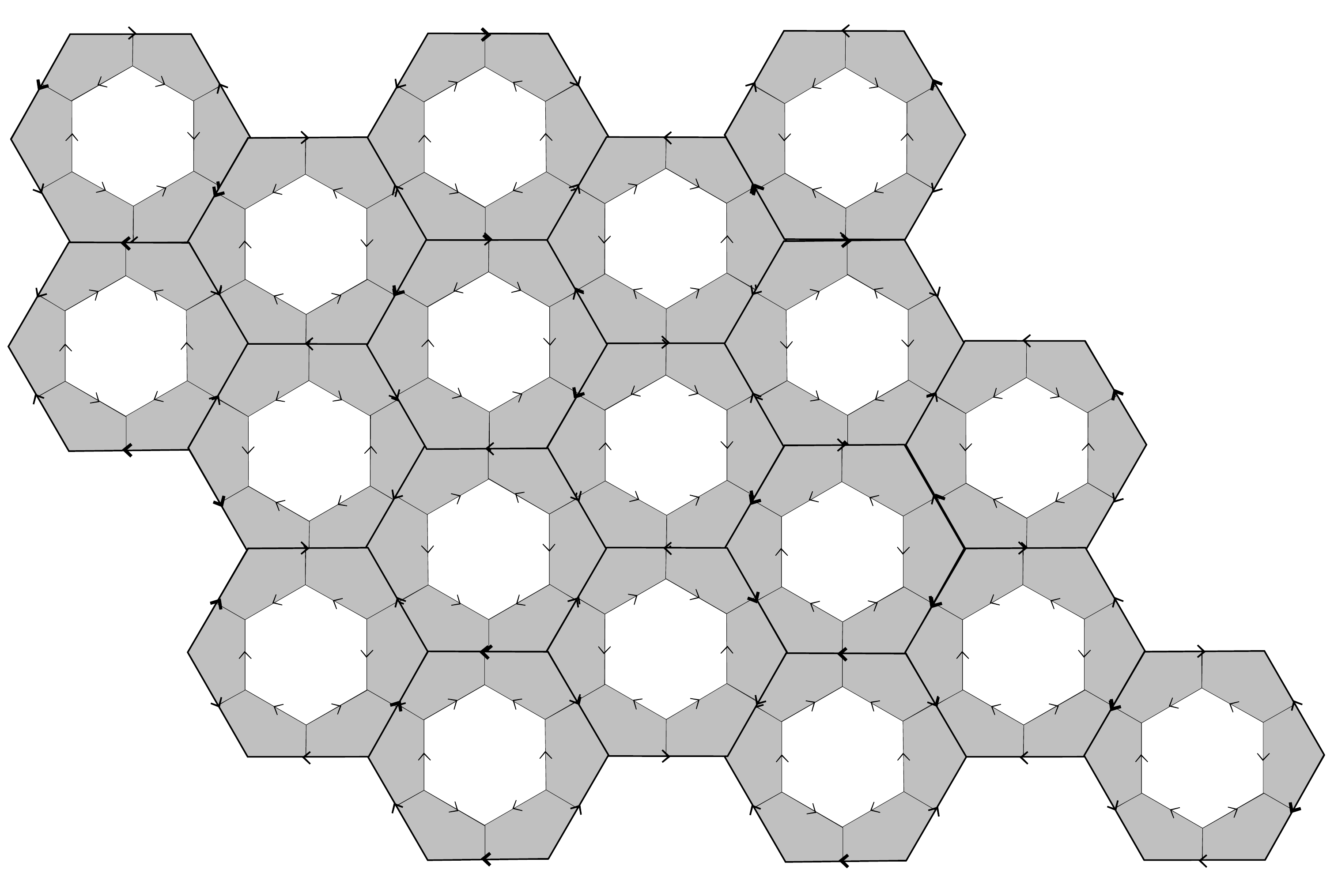}
\caption{A legal patch of double hexagon tiles. Where three double hexagons meet at a vertex, a corner hexagon is
created. Note the mixed arrowing of these small grey corner hexagons. }
\label{doubleHexPatch}
\end{figure}

Suppose that all the corner hexagons of the tiling can be well-arrowed in this way. Collectively the inner hexagons together with the corner hexagons form another hexagonal tiling of the plane, if we ignore
the question of their arrows matching. However, there does arise the question of whether or not all the common edges of adjacent
small hexagons actually do have matching arrows, that is, whether or not this new tiling is well-matched. The double hexagon tiling is called {\bf legal} if they do. Thus
a double hexagon tiling is {\bf legal}
if its outer hexagons are well-matched, all of its corner hexagons can be completed to well-arrowed hexagons, and 
the consequent well-arrowing of the small hexagons completes 
to a small hexagon tiling of the plane which is well-matched.
In this situation we have two well-matched hexagonal
tilings, one using the large hexagons and the other using the small ones inner and corner hexagons.
Fig.~\ref{doubleHexPatch} shows a patch of a legal double hexagon tiling. 

The double hexagon tiles that we are discussing are also called {\bf Penrose hexagonal tiles}, and a legal tiling is a {\bf Penrose hexagonal tiling}. We shall use both names in this paper, because within the 
context of understanding the intimate relationship between Penrose hexagonal tilings (based on the large hexagons) and Taylor-Socolar hexagonal tilings (based on the
small hexagons), it is convenient at times to simply think in terms of legal double hexagon tilings. 

%%%%%%%%%%%%%%%%%%%%%%%%%%%%%%%%
\section{Decorations and triangles}
\begin{figure}[h]
\centering
\includegraphics[width=0.6\textwidth]{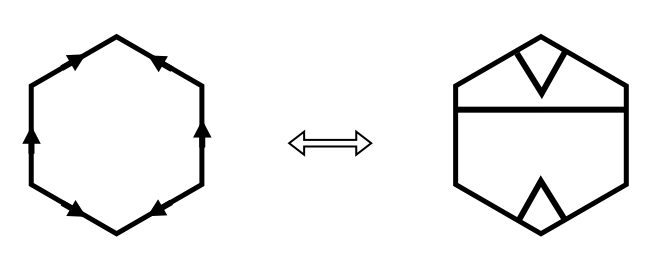}
\caption{ Well-arrowed hexagonal tiles can be converted into hexagonal tiles with stripes. These 
decorations fit together to make triangulations of the plane.}
\label{basicArrowedHexagon}
\end{figure}

%%%%%%%%%%%%%%%%%%%%%%%%%%%%%%%%
There are other decorations of well-arrowed hexagons and double hexagon tiles that are equivalent representations
of the arrowing but help to make the underlying geometry of the tilings more transparent. 
The first of these is the marking
of well-arrowed hexagons shown in Fig.~\ref{basicArrowedHexagon}, which replaces the arrows of a well-arrowed hexagon with a black stripe and two black corner markings. Initially we will use this representation of the arrowing with the small 
hexagons, later for the outer hexagons. 

Notice that when two well-arrowed hexagons are attached along some edge so that the corresponding
arrows match (i.e. they are well-matched), the black markings line up, either stripe to stripe, stripe to corner,  or corner to corner, to create
an extended black path. In fact, that the stripes and corner markings match to form extended paths, is exactly
the same as arrow matching. In the resulting paths
the corner markings indeed serve as corners at which the direction of the path changes. If we have a well-matched tiling of
well-arrowed tiles, we will have also a set of paths. It is easy to see that if, purple in following a path, it turns right or left at a corner,
then at its next corner it will turn in the same sense (again right or again left) and so the resulting paths 
will be equilateral triangles (the corners create $60^{\circ}$ angles). 

The only way this can fail is if there are paths that extend infinitely in some direction
along some straight line. Such a tiling is called a {\bf singular} (or {\bf non-generic}) tiling. Later on we will examine the similar
paths created by the stripes and corner markings on the large hexagon tiles, and the same issue of
singularity will arise.  

The generic situation is that of {\bf non-singular} (or {\bf generic}) tilings,
that is, all of the paths form triangles. 
In this paper, in order to keep all the essential ideas clear, we shall always assume non-singularity, though at this point we need it only with the small hexagons and their markings.  The resulting  triangles come in various sizes and arrangements, and this is something we address in the next section.

\begin{figure}[h]
\centering
\includegraphics[width=0.8\textwidth]{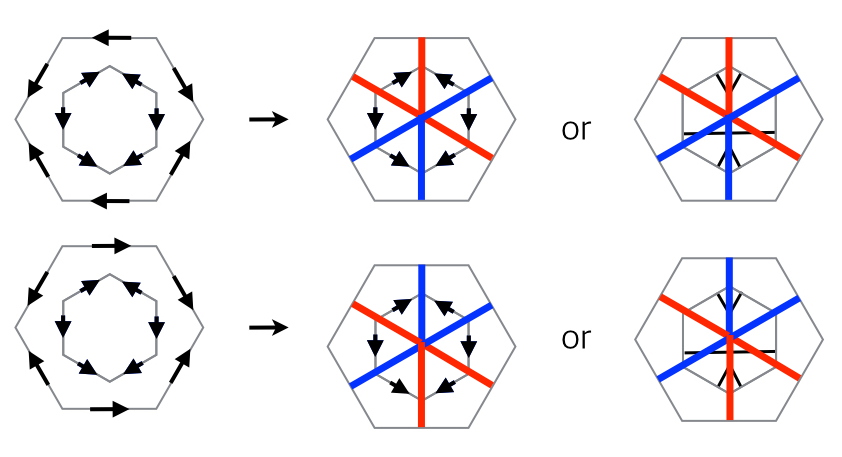}
\caption{Arrows on the edges of the outer hexagon which are oriented in the counterclockwise direction are represented by short red half-diameters. For arrows in the clockwise orientation we use short blue half-diameters. If we include the striping decoration
of the inner hexagons as well (Fig.~\ref{basicArrowedHexagon}), we arrive at the fully decorated double hexagons shown on the right-hand side of the figure. Evidently the decorated tiles carry information fully equivalent to the arrowing. Notice
that proper matching of the edges of outer hexagons is equivalent to a {\bf change} of color as the short diagonal passes
through the common edge. Also notice that if one holds the black stripe horizontally, then as one moves along a full blue diagonal from right to left, the diagonal passes through the black stripe from above to below. It is the other way around for red stripes.  We use this observation to make the color determination of the short diagonals
in Fig.~\ref{bigHexLineUp}.}
\label{NewDeco-on-PenroseTile}
\end{figure}

The second decoration is one that we make to double hexagon tiles, replacing the outer arrowing by colored short
diagonals (short diagonals are the ones that pass at right-angles between opposite edges, as opposed to long diagonals
that pass from vertex to opposite vertex). This is explained in Fig.~\ref{NewDeco-on-PenroseTile}.

If we begin with a legal double hexagon tiling then we know that we end up with two well-matched hexagonal tilings:
one of large hexagons and one of small hexagons.  Since we are assuming non-singularity, the well-matching of the small hexagons leads to a collection of triangles
on the plane -- equilateral triangles created by the stripes on the small hexagons. Each small hexagon has an inner part of some edge of a triangle across it and the corners of two other triangles, one on each side of that edge, so altogether
each small hexagon is involved with three triangles. 

The very smallest triangles (called level $0$ triangles) are those composed by putting three
corner markings together around a common vertex of three small hexagons. Every stripe in a hexagon obviously belongs to a triangle
larger than these smallest ones. Indeed there are triangles of ever increasing sizes, without limit. It is this result that
we will establish in the next section.

%%%%%%%%%%%%%%%%%%%%%%%%%%%%%%%%
\section{Nesting and hierarchy} \label{Nesting-and-hierarchy}
%%%%%%%%%%%%%%%%%%%%%%%%%%%%%%%%

Let us continue with a non-singular legal double hexagon tiling $\cD$, in which we have completed its small hexagons to a well-matched hexagonal tiling and then resolved everything into triangles by decorating each of the small hexagons. 

\begin{figure}[h]
\centering
\includegraphics[width=0.8\textwidth]{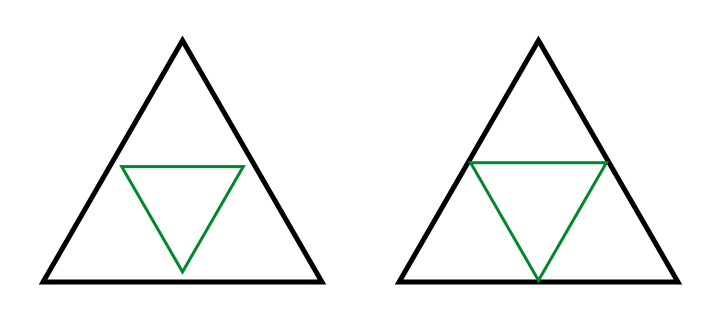}
\caption{The green (smaller) triangle is nested in the larger (black) one. Notice the two ways of drawing this.
A patch of triangles that arise from our tiling presents the inner triangles as being totally in the interior of the outer ones, as shown on the left-hand side here. We
often will wish to allow the inner triangles to stretch to meet the boundaries of the outer triangles, as shown
on the right hand side of the figure. This creates three new triangles called {bf corner triangles}, so that the outer triangle is now decomposed as four equal sized triangles.}
\label{nestedTriangleScheme}
\end{figure}

A triangle is {\bf nested} in another one if it appears as in Fig.~\ref{nestedTriangleScheme}, where
the smaller triangle is nested in the larger. In this section we will prove that except for the very
smallest triangles (the ones made from three corners) every triangle has another one nested inside it.
From this, we will see that every triangle has inside it a sequence of triangles nested within each other, diminishing
in size to the smallest size triangles. We refer to this phenomenon by saying that all triangles are {\bf nested within}.

There is more to this. Let us stretch out, or expand, the 
triangles created by the decorations of the small hexagons (inner hexagons and corner hexagons) so that the corners meet the edges of the triangle surrounding them as illustrated in Fig.~\ref{nestedTriangleScheme}. In doing this each nested triangle produces three neighbors that
fill out the whole triangle that it lies in. 
In fact there is a nesting that involves one triangle sitting inside another of exactly twice the linear size, so that the larger triangle is decomposed into four equal-sized equilateral triangles of which the nested triangle is one. The three triangles that emerge as neighbors of the nested one are called {\bf corner triangles} (not to be confused with
the smallest triangles that we formed out of three corner markings).

\begin{figure}[h]
\centering
\includegraphics[width=0.7\textwidth]{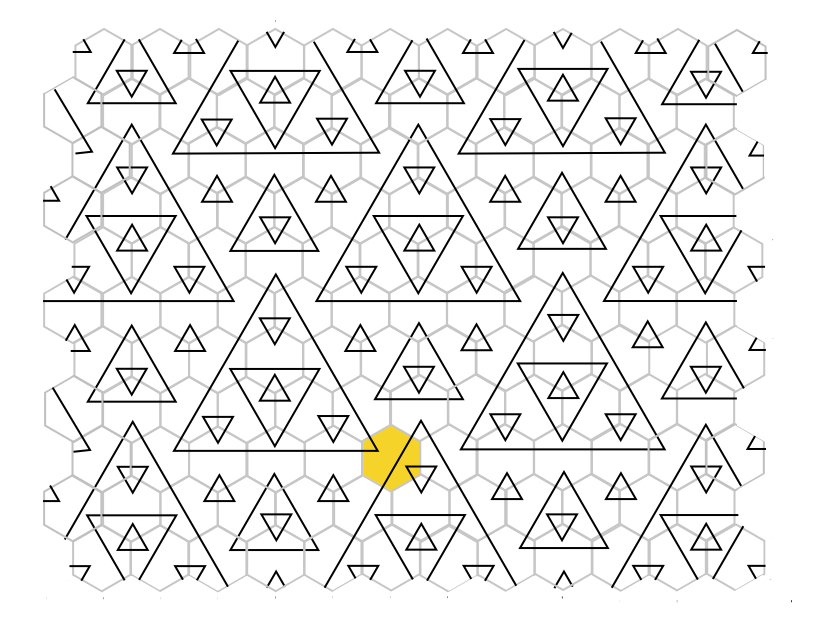}
\caption{Periodic arrangements of triangles are possible if only decorations of the small hexagons are
used. }
\label{TrianglePatch-008}
\end{figure}

We will speak of the patterns of 
triangles (expanded or not) which are formed by the decorations of the small hexagons as {\bf arrangements of triangles} and derive their nesting properties as we proceed.
Notice that without the implications derived from the decorations of the outer hexagons, it is possible to get an arrangement of triangles like the one shown in Fig.~\ref{TrianglePatch-008}, which is visibly periodic.

In all we shall see that triangles that are nested within appear on ever increasing scales, so
there is a hierarchical structure.  We shall call such an arrangement of triangles a {\bf nested triangulation}.

When we refer to the sizes of triangles in 
one of arrangements of triangles, we will always refer to the side lengths of the 
stretched out versions. 
Lengths are normalized so that the smallest triangles will be 
of side length equal to $1$. We will see that with this normalization all lengths are powers of $2$.

In the sequel we will commonly use both versions of the triangles and nested triangles that emerge out of
 our discussion--the original ones that come
from the decorated hexagons, and the stretched out ones that give us 
the arrangements of triangles.
Once we have proved that all triangles are nested within, it is trivial to convert from one picture to the other.

\begin{figure}[h]
\centering
\includegraphics[width=0.8\textwidth]{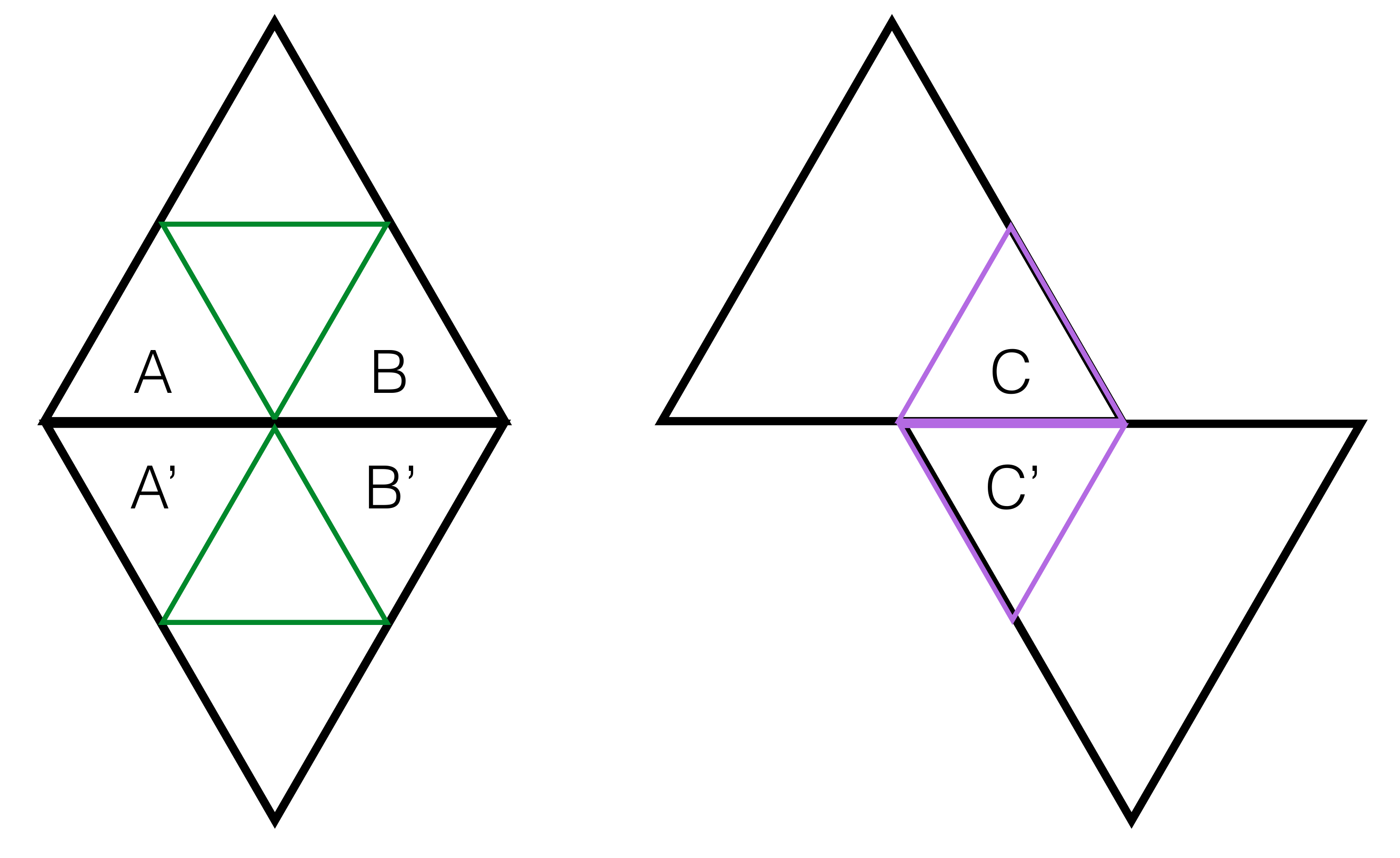}
\caption{On the left side the black triangles form an opposite pair of triangles. So too do the pairs
of smaller triangles labelled $A,A'$ and $B,B'$. On the right-hand
side the two black triangles do not form an opposite pair, but the two purple triangles $C,C'$ do. }
\label{oppositePairs}
\end{figure}

In the stretched out version, two triangles are said to make an {\bf opposite pair} if they are of the same size and share a common
edge, see Fig.~\ref{oppositePairs}. Notice that there is no specification of how each of the two opposing triangles fits into the overall 
arrangement of triangles.

\begin{prop}\label{nesting}
Let $\cD$ be a non-singular legal double hexagon tiling. Complete its small hexagons to a well-matched
hexagonal tiling and let $\cT$ be the resulting arrangement of expanded triangles formed from the decorations of the small hexagons. Then \\
{\rm{(i)}} all triangles occur in opposite pairs;\\
{\rm{(ii)}} the side lengths of the triangles are all of the form $2^{k}$ for some $k=0,1,2,\dots$ 
{\rm (}$k$ is called the 
{\bf level} of the triangle\rm{)};\\
{\rm{(iii)}}  every triangle is nested within.
\end{prop}

The proof of Prop.~\ref{nesting} is by induction on the size of triangles. 
The smallest triangles have side length $1= 2^{0}$ (level $0$). There is no
nesting within to take place. The stretched triangle pattern created by these triangles is in itself 
a genuine triangular lattice of the plane and, in particular, every triangle edge borders an opposite pair of triangles. 

We now assume that the three statements of Prop.~\ref{nesting} have been proved up to 
some level $k$. 

\begin{figure}\centering
\includegraphics[scale=0.4]{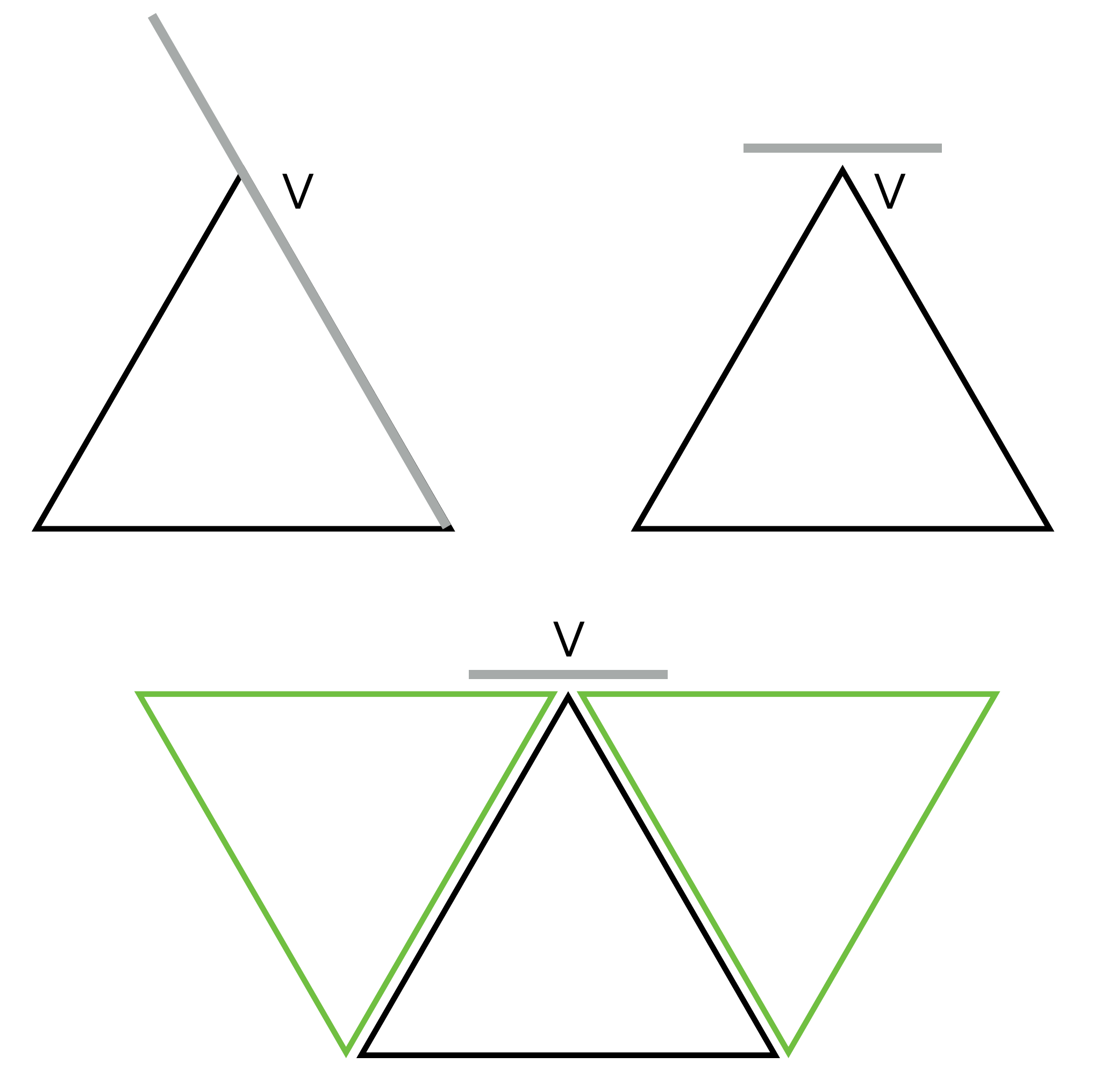}
\caption{$v$ is a vertex of a triangle (shown in black) of level $k$. At $v$ there is a small hexagon, and its stripe allows
for only two things to happen: either one of the edges of the triangle extends through $v$, in which case it
is an edge from a larger triangle, or there is an edge passing through $v$ that is parallel to the opposite
edge of the triangle. In the latter case we use (i) to place down the two opposite pairs of triangles  
of adjoining triangles, shown in green (the adjacent edges are actually coincident edges of course). Then
we see that the edge through $v$ must actually be an edge that includes both the top edges of the
green triangles: thus again a larger triangle.}
\label{largerTriangles}
\end{figure}

First we check that there must be triangles of size larger than $2^{k}$. Fig.~\ref{largerTriangles} shows why. It shows
a triangle of level $2^{k}$ and uses matching triangles to see that there must be larger ones. 

\begin{figure}\centering
\includegraphics[scale=0.5]{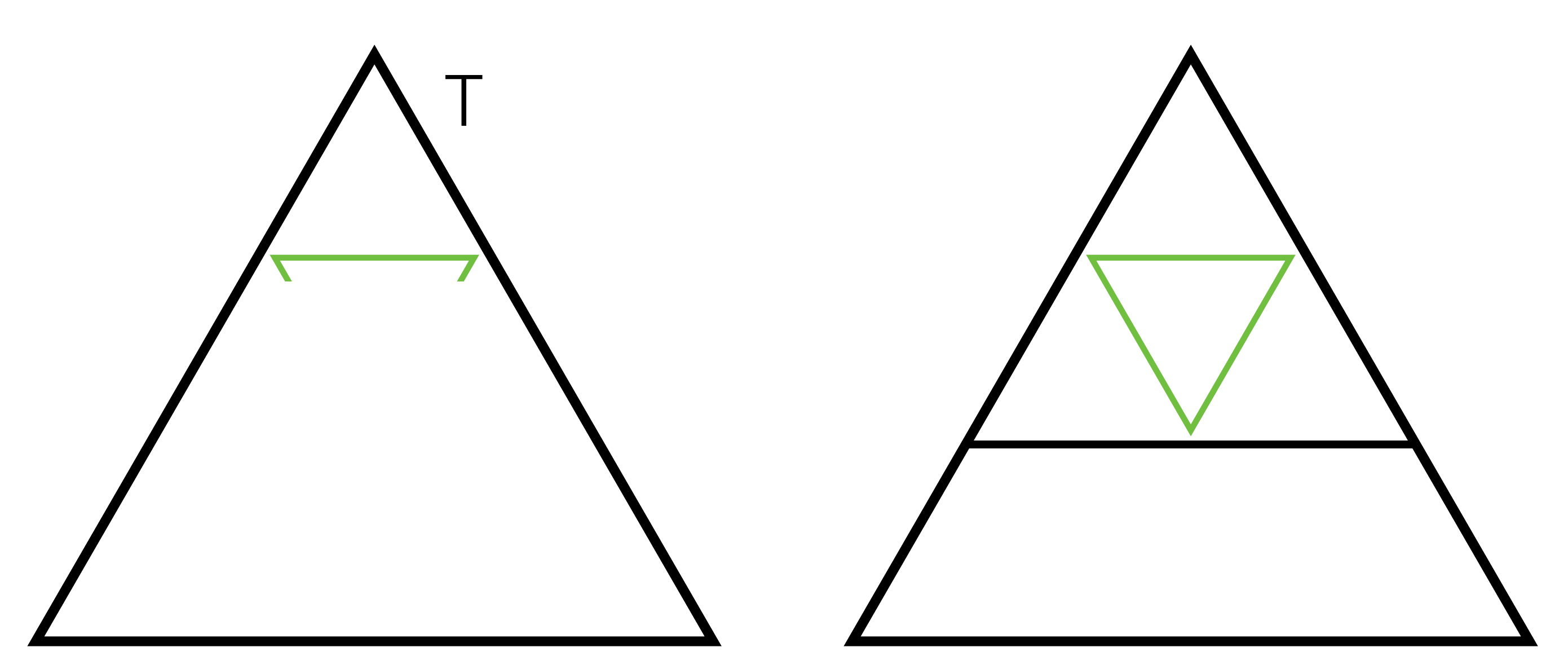}
\caption{We are given a triangle $T$ of minimum size larger than $2^{k}$. From one of its vertices
(we have taken it as the top one here) we fit in the largest sub-triangle possible (at the very least, there is always a triangle of level $0$ that can be fitted in). Its lower side is indicated in green. It must turn inwards at the sides of
$T$ and complete to the opposite green triangle. By the induction assumption its other two sides
also complete to opposite pairs, and this leads to the new black triangle with the green triangle nested in it. 
Since we started from a maximal sized sub-triangle, this larger black triangle must in fact the entirety
of our original triangle $T$. This shows that $T$ has edge length $2^{k+1}$. The visible nesting and the induction hypotheses
show that the new triangle nested within.}
\label{internalNesting}
\end{figure}

We now take any triangle $T$ of the next size, say $m$, that is larger than $ 2^{k}$. 
We see immediately that $m= 2^{k+1}$ and it is internally nested, in the right way, Fig.~\ref{internalNesting}.
This completes the induction steps for parts (ii) and (iii).

\begin{figure}\centering
\includegraphics[scale=0.4]{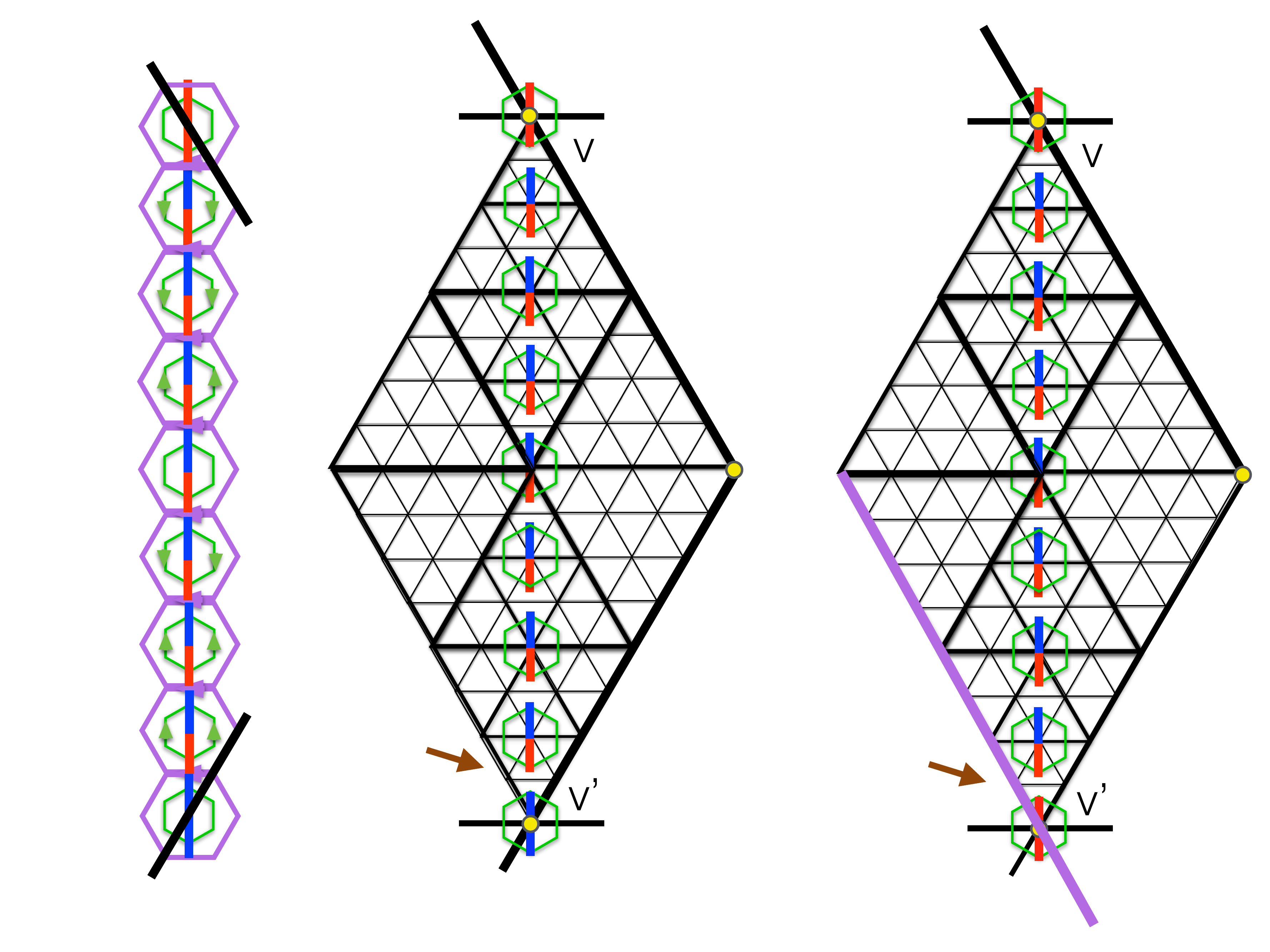}
\caption{In the center of the figure we see the large pair of triangles making a diamond shape between the
two extreme vertices $v$ and $v'$, which are assumed to be vertices arising from centers of
the large hexagons of the double hexagon tiling. The diamond is made up of an opposite pair of fully internally nested 
triangles. The triangles are all in their stretched form, but the line thicknesses indicate
the nesting relationships. On the left side we see the corresponding arrangement of double hexagons that surround
the small hexagons between $v$ and $v'$. The arrows must
match, but their common orientation in the horizontal direction is irrelevant here. There is a color change as
we cross each triangle edge at right angles. The key point is what happens at the ends of the diamond
as the color line crosses edges (indicated by the thickest lines) which are not at right angles to it. The main edge at $v$ is shown by
the heavy black line. The rules for coloring hexagons show that the color stripe is fully red here, see 
Fig.~\ref{NewDeco-on-PenroseTile}. The main edge at $v$ has to be
matched with its partner at $v'$, where color strip changes to fully blue.  Notice the correct color change
at the arrow. The scenario shown at the right side of the Figure, where the main edge at $v'$ is shown in purple, cannot occur because of the color change violation at the brown arrow. }
\label{bigHexLineUp}
\end{figure}

We now come to the proof of part (i). It is useful to prepare this by looking at the situation pictured in
Fig.~\ref{bigHexLineUp}. What this shows is how the coloring of the tile decorations is related to the matching of opposite
triangles. The color rules show that as a color diagonal crosses a triangle edge at right-angles it changes color.
When it crosses an edge that is not at right angles to it then it does not change color, but, as
we have noted in the caption to Fig.~\ref{NewDeco-on-PenroseTile}, its color is related to the way in which it crosses the edge. This figure is the basis for our proof of matching triangles.

\begin{figure}\centering
\includegraphics[scale=0.5]{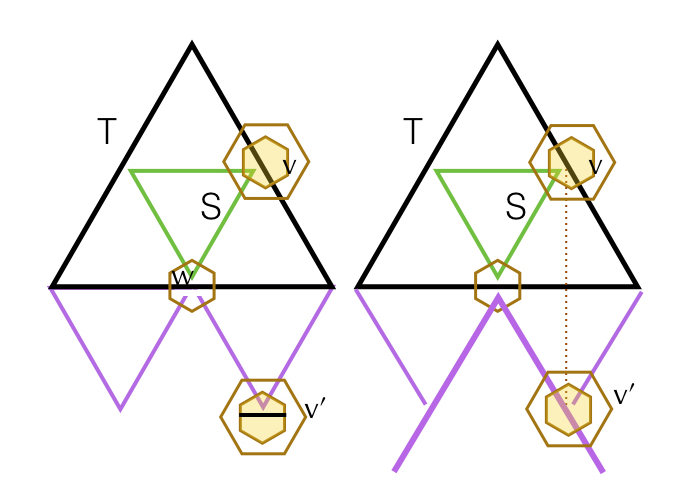}
\caption{The vertex $v$ is the unique vertex of triangle $S$ which is the center of a double hexagon tile. The two
lower corner triangles formed in $T$, with common vertex $w$, have opposite triangles shown in purple. The question
is, what happens at the vertex $v'$? Neither of the two possibilities shown here can occur. On the left, the triangle
with vertices $v'$ and $w$ would be left unclosed at $w$. On the right we are in the situation shown on the right side of
Fig.~\ref{bigHexLineUp}, which we know violates the color change property.} 
\label{twoGoodSides-I}
\end{figure}

Continuing to the proof of (i) we start with a triangle $T$ of level $k+1$ and show that it must be matched 
by a triangle of the same level on each of its sides. Let $S$ be the largest equilateral triangle nested
in it. Now, any equilateral triangle of any level $2^{r}$ in our arrangement of triangles
has exactly one vertex at the center of a large hexagon of the double hexagon tiling. This has to do with 
the lattice structure induced by the arrangement of triangles coming from double hexagon tiling, and though pretty self-evident in the figures,
is explained algebraically in \S\ref{algNt}. In Fig.~\ref{twoGoodSides-I}
we have made such a choice, indicating it by the small yellow circle at $v$. We shall use this coloring convention to mark other
vertices that are centers of the large hexagons. We shall start by showing that there must be matching
triangles to $T$ on the two sides of $T$ on which $v$ does not lie. 

\begin{figure}\centering
\includegraphics[scale=0.5]{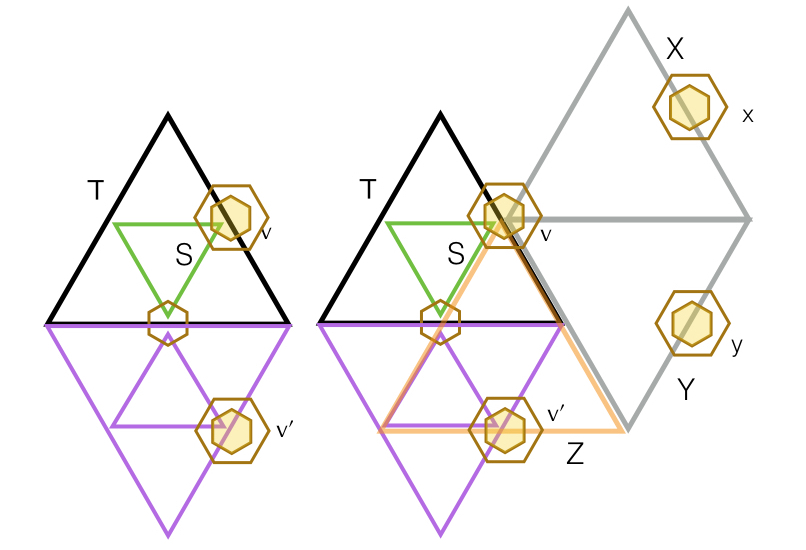}
\caption{The left-hand side shows the matching that has to take place on the lower side of $T$.
 The right-hand side shows what happens if there is not an opposite match to $T$ on the
 side containing $v$. This corresponds to the right-hand side of 
 Fig.~\ref{twoGoodSides-I}. Chasing around the pairs of opposite matching triangles yields $X,Y,Z$ and the latter
 is clearly totally mis-matched with $T$. }
\label{thirdGoodTriangle-I}
\end{figure}

The triangle $S$ creates a partition of $T$ into itself and three surrounding triangles, and we know that each of these
must have an opposite match. We show these matching triangles along the lower edge of $T$ solid edge. The shape of the small hexagon
at their intersection $w$ must be of the type shown in the Fig.~\ref{twoGoodSides-I}. What does the small hexagon at $v'$ look like? The caption to 
Fig.~\ref{twoGoodSides-I} shows that neither of the two possibilities shown there is possible.  
Thus the remaining possibility, which is that of a matched triangle for $T$, must occur. This is illustrated on the left side
of Fig.~\ref{thirdGoodTriangle-I}. 
This same argument can be applied to the other side of $T$ which does not contain the point $v$. 

There remains the task of proving that the side that contains $v$ also matches $T$ to a triangle of the same size. Let us suppose
that on this side the matching fails. The argument we have just used tells us that in this case on this side
we will see a triangle $X$ which aligns its corner at the point $v$. The point $x$ is the center of a double hexagon, just
like $v$ was, so it follows by what we have proved that the triangle $Y$ shown exists and matches it. Again it has
a point $y$ which is a double hexagon center and so $Y$ produces the matching triangle $Z$. But this is clearly a contradiction
since $Z$ overlaps but does not coincide with the triangle $T$. This contradiction shows that there is a matching
opposing triangle along the edge of $T$ containing $v$.

With this we conclude part (i) of Prop.~\ref{nesting}, and so the entire proposition. \qed

\bigskip
Since there are triangles of every level, it is impossible that there are any translational symmetries.

\begin{prop}\label{aperiodicity}
Every non-singular legal double hexagon tiling is aperiodic. \qed
\end{prop}

\begin{figure}\centering
\includegraphics[scale=0.5]{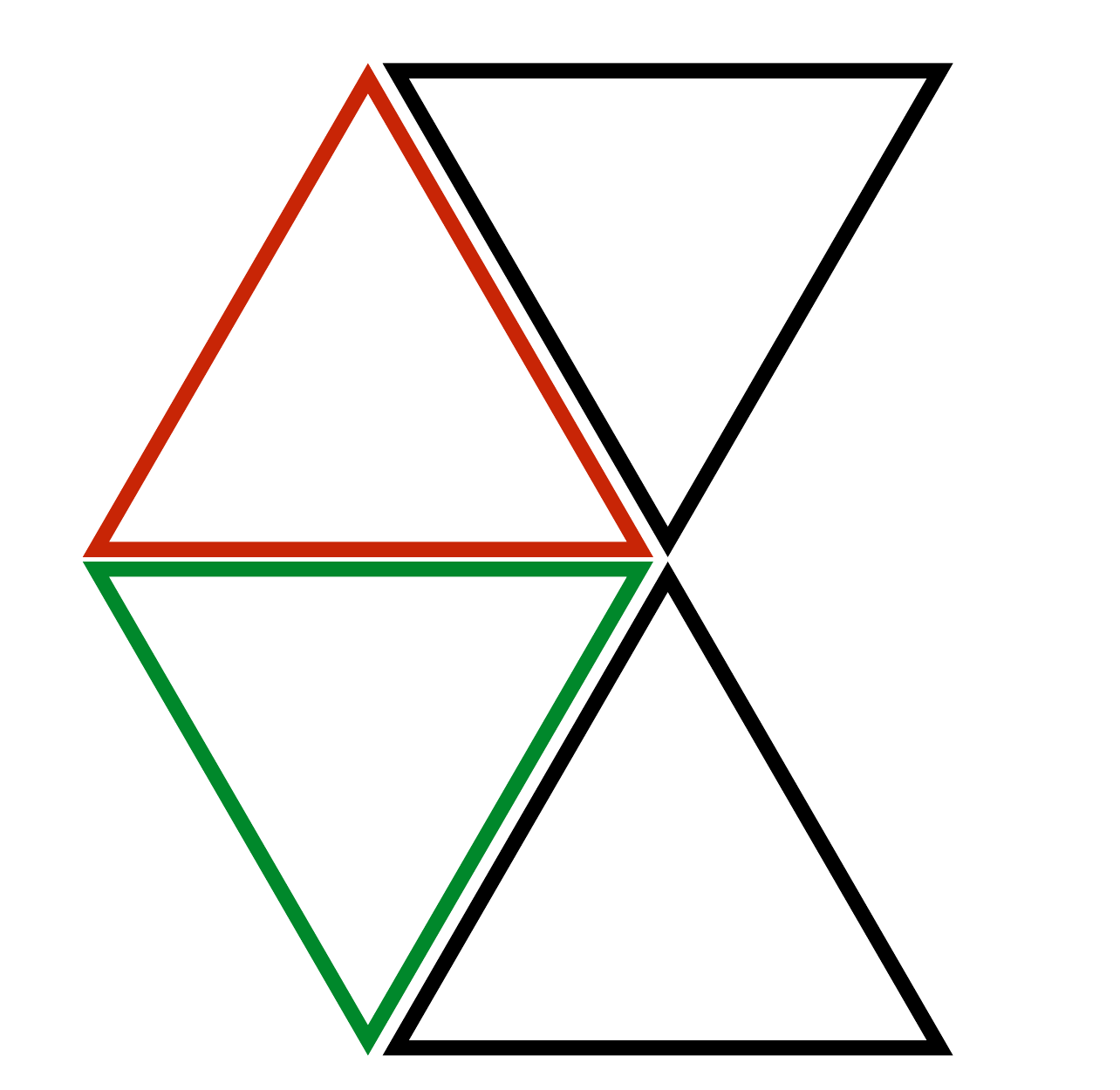}
\caption{The two triangles whose corners make up the pair of corner markings on a well-arrowed
tile are always of the same size. Here the pair of triangles is shown in black. The matching on opposite
sides of triangles leads to the red and green matched triangles. Since these too must match, we see
that all four triangles are of the same size. Note that there is no presumption here about how these triangles
lie in larger triangles.}
\label{matchingCornerTriangles}
\end{figure}

Looking at Fig.~\ref{matchingCornerTriangles} we see:

\begin{prop} \label{matchingCorners}
In any small hexagonal tile the triangles that arise from its two opposite corner markings are
of the same size. \qed
\end{prop}

For future reference we also note (see Fig.~\ref{nestedTriangleScheme}):
\begin{prop}\label{centralEdgeProperty}
In the nested triangulation created by the standard edge and corner markings of 
arrowed hexagonal tiles, every stripe forms the central part of an edge of a some triangle.
\end{prop}

%%%%%%%%%%%%%%%%%%%%%%%
\section{The algebra of nested triangulations}\label{algNt}
%%%%%%%%%%%%%%%%%%%%%%%

\begin{figure}
\centering
\includegraphics[scale=0.4]{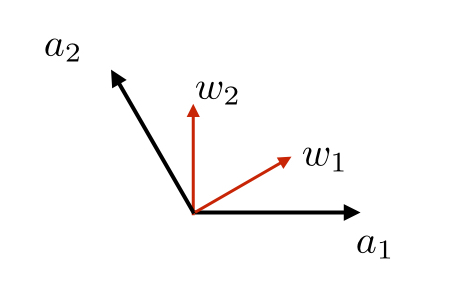}
\caption{The basis vectors for the lattices $Q$ and $P$. The directions of edges in the $Q$-triangulations
are $\pm a_{1}, \pm a_{2}, \pm(a_{1}+ a_{2})$, which are called a-directions, and those of $P$-triangulations are
$\pm w_{1}, \pm w_{2}, \pm(w_{2}- w_{1})$, which are called w-directions.}
\label{aAndW-I}
\end{figure}

\begin{figure}
\centering
\includegraphics[scale=0.3]{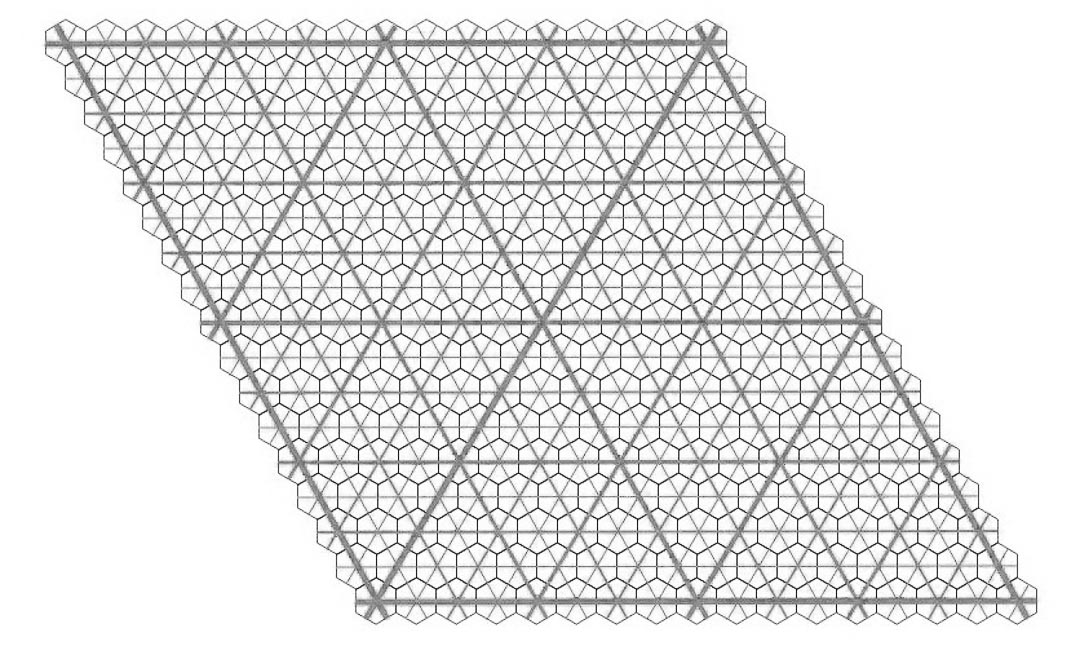}
\caption{A nested triangulation of the plane. Triangulations of increasing scales (four are shown here)
coexist in one underlying triangulation. Each increase in scale can be created in four different ways 
by choosing one vertex of the previous scale as a reference point.}
\label{basicTriangulation}
\end{figure}

If we start with a non-singular legal double hexagon tiling then we obtain a tiling of the plane with the small
hexagons. The centers of these hexagons form a triangular lattice of the plane composed of level $0$ (side length $1$) equilateral triangles, as we have seen. For definiteness
we now specify this lattice as a set of points in $\R^{2}$, namely the set of points
$Q= \Z a_{1}+ \Z a_{2} \subset \R^{2}$, where
$a_{1}= (1,0)$, $a_{2}= (-\frac{1}{2}, \frac{\sqrt{3}}{2})$, Fig.~\ref{aAndW-I}. Joining nearest neighbors of $Q$ produces the triangular lattice of level $0$ triangles, indicated by the thin lines in Fig.~\ref{basicTriangulation}.

We know that there are also triangles of level $1$ (side length $2$). They are matched across each of their edges, and so there
is a second triangular lattice of the plane by equilateral triangles. This meshes precisely with the first, in the sense
that each level $1$ triangle is composed of four level $0$ triangles. The vertices of the level $1$ triangles
form a coset $q_{1} + 2Q$ of $Q$. 

We may repeat this process, now looking at triangles of level $2$, whose vertices lie on a coset 
$q_{1}+ q_{2} + 4Q$ (where $q_{1}\in Q$ and $q_{2} \in 2Q$). Continuing this way we led to view our
nested triangulations in terms of ever refined cosets from the sequence 
\[Q \supset 2Q \supset 4Q \supset 8Q \supset \cdots \,.\]
Thus the double hexagon tiling leads to the sequence
\[{\bf q} = (q_{1}, q_{1}+ q_{2}, q_{1}+ q_{2}+ q_{3}, \dots \ ) \, ,\]
where each $q_{k}\in 2^{k-1}Q $. We refer to such a sequence as 
a {\bf $Q$-nested triangulation} $\cT$. Indeed we see that any such sequence corresponds to a sequence
of triangular lattices with each level nested within the next, that is to say every triangle of one
level appears as a corner triangle or as a central triangle within a triangle of one level higher.
Specifically, up to rotation a typical triangle of level $k$ has vertices $x, x+2^{k}a_{1}, x+2^{k}a_{2}$,
all of them lying in one coset $q_{1}+ \cdots + q_{k} +2^{k}Q$, where we assume $k\ge 1$. The
mid-points of its edges are $x+2^{k-1}a_{1}, x+2^{k-1}a_{2}, x+2^{k-1}a_{1}+ 2^{k-1}a_{2}$, which form the
vertices of a triangle
of level $k-1$ in the coset $q_{1}+q_{2}+ \cdots + q_{k-1} +2^{k-1}Q$, hence the nesting.

The sequence $\bf q$ can be interpreted as an $Q$-adic element of the
inverse sequence
\begin{equation}\label{Qzero}
 \overline {\bf Q} :Q/Q \leftarrow  Q/2Q \leftarrow Q/4Q \leftarrow Q/8Q \leftarrow \cdots \,.
 \end{equation}
This sets up a bijective correspondence between $Q$-adic numbers and nested triangulations,
and we will write $\cT = \cT(\bf q)$ when we wish to make the connection explicit. In the sequel 
$\overline {\bf Q}$ will be written as $\overline {\bf Q_{0}}$, the first in a series of such inverse limits. 

The condition of non-singularity has  algebraic consequences. The $Q$-nested
triangulation is singular
if some of the paths created by the stripes of the small hexagons do not close, but rather extend indefinitely.
Since the directions of the stripes are all in the $a$-directions of the lattice $Q$ (see Fig.~\ref{aAndW-I}) this is equivalent to saying that there is an infinite path of 
edges consecutive edges in some direction $a \in \{\pm a_{1},\pm a_{2}, \pm(a_{1}+ a_{2})\}$
and this in turn implies that ${\bf q}$ lies in $x + \overline{\Z_{2}} a$ for some $x\in Q$. Here 
$\overline{\Z_{2}}$ is the $2$-adic integers. We need to avoid 
${\bf q}$ having this form. See \cite{TSMS} for more details. 

In order to interpret the double hexagon tiles in this algebraic setting, we need, along with $Q$, its $\Z$-dual $P$, relative to the standard dot product on $\R^{2}$.
Thus $P= \Z w_{1}+ \Z w_{2}$ where $w_{1}=\frac{2}{3}a_{1}  +\frac{1}{3}a_{2} $ and
$w_{2}=\frac{1}{3}a_{1}  +\frac{2}{3}a_{2} $, Fig.~\ref{aAndW-I}. We note that
\[ P\supset Q \supset 3P \supset 3Q\supset 9P \supset \cdots \,,\]
all the steps being of index $3$.
 In fact each of the lattices in this chain is a scaled and rotated version of the one before it,
 and in particular a scaled and rotated version
 of the original triangular lattice $Q$. They are all triangular lattices.
 We do not use $P$ directly in what follows, but rather $3P$, since
 we wish to keep everything inside the initial lattice $Q$.

There are two clear differences between the triangular lattices arising from $Q$ and $3P$. The first is
that the basic triangles in $3P$ have side length $\sqrt 3$, so factors of $\sqrt 3$ relate scales of
$Q$- and $3P$-nested triangulations. The second is that the directions of the sides in all the $Q$-nested triangulations
(at all scales) are $\{\pm a_{1}, \pm a_{2}, \pm (a_{1}+a_{2})\}$, which we are referring to as
$a$-directions, and those for $3P$-nested triangulations are $\{\pm 3w_{1}, \pm 3w_{2}, \pm 3(w_{2}-w_{1})\}$, which we
have called $w$-directions. These two sets of directions are interchanged by $90^{\circ}$ rotations.

In a legal double hexagon tiling the large hexagons share one third of the vertices of the small hexagons,
and their centers form some coset $c+3P$ of $Q \mod 3P$. This brings us to a second inverse sequence of groups
based on $3P, 6P, 12P, \dots$ and corresponding group $\overline{\bf Q_{1}}$ which is related to $\overline {\bf Q}$
as shown in the commutative diagram \eqref{commutativeDiagramSmall}. All the mappings here are the natural homomorphisms that arise from factoring out larger subgroups.

 \begin{align}\label{commutativeDiagramSmall}
\overline{\bf Q_{0}}: &\quad &Q/Q  & \quad \longleftarrow & Q/2Q & \quad \longleftarrow & Q/4Q & \quad \longleftarrow &   Q/8Q & \leftarrow\cdots  \nonumber\\
\uparrow &&\uparrow & & \uparrow & & \uparrow & &  \uparrow & &\nonumber\\
\overline {\bf Q_{1}}: & \quad &Q/3P  & \quad \longleftarrow & Q/6P & \quad \longleftarrow& Q/12P & \quad
\longleftarrow & Q/24P &\longleftarrow \cdots \nonumber \\
\end{align}

Given the choice of $\bf q \in \overline Q$ and $c \in Q/3P$, there is a unique element 
\[{\bf r} = (r_{1}, r_{1}+ r_{2}, r_{1}+ r_{2}+ r_{3}, \dots \ ) \, ,\] 
where $r_{k} \in 2^{k-1}3P$ for each $k$,
which maps onto $\bf q$ and has $r_{1} \equiv c \mod 3P$. This follows from the more general fact:

\begin{lemma} \label{subgroupIntersections} For all $k,l \in \N$,
\[ 3^{k}2^{l-1}P \, \cap \, 3^{k-1}2^{l}Q = 3^{k} 2^{l} P \quad \mbox{and}\quad
 3^{k}2^{l}P \, \cap \, 3^{k}2^{l-1}Q = 3^{k} 2^{l}Q
\, .\]
\end{lemma}

\noindent Proof: Dividing out common powers of $2$ and $3$, we are reduced to proving that
$3P \cap 2Q = 6P$
and $2P\cap Q= 2Q$, respectively, both of which are trivial to verify.
\qed

\medskip

So, given a non-singular legal double hexagon tiling, we arrive not only at an element
$\bf q \in \overline{Q_{0}}$ but also an element $\bf r \in \overline{Q_{1}}$. This element picks out
one coset $r_{k} + 2^{k}3P$ for each $k=0,1,2, \dots$ and so should determine a family $\cT^{+}(\bf r)$ of nested
triangulations, just the same way as $\bf q$ did. To guarantee that we really do have a non-singular $3P$-nested triangulation, that is
to avoid infinite lines, we have to make an assumption similar to the non-singularity assumption 
that we have already seen. This time the triangle edges are in $w$-directions, so we must assume
that ${\bf r}$ does not lie in $x+ \overline{\Z_{2}}3w$ for any $w\in \{w_{1}, w_{2}, w_{2}-w_{1}\}$.

Thus the joint conditions equivalent to non-singularity are that for all $x\in Q$,
\begin{itemize}
\item[(i)] ${\bf q}$ does not in $x + \overline{\Z_{2}} a$ for any $a \in \{\pm a_{1},\pm a_{2}, \pm(a_{1}+ a_{2})\}$ \,;
\item[(ii)] ${\bf r}$ does not lie in $x+ \overline{\Z_{2}}3w$ for any $w\in \{w_{1}, w_{2}, w_{2}-w_{1}\}$ \,.
\end{itemize}
These are the same conditions as appeared in \cite{TSMS}, though we did not use
double hexagon tilings there.
We will pursue the detailed study of the singular double hexagon tilings in another paper.

\medskip
Returning to our discussion, the situation is this. We are given a generic double hexagon tiling whose small hexagons are centered on $Q$ and whose large hexagons are centered on $c + 3P$. We thus have two nested triangulations $\cT({\bf q})$ and
$\cT^{+}({\bf r})$, the first being determined by the markings on the small hexagon tiles
and the second simply by the algebra of the commutative diagram \eqref{commutativeDiagramSmall}.
 Since the large hexagon tiles can be given their own stripe and corner markings based on their arrowing in just the same way as we did for the smaller hexagons, it is natural to
ask whether or not this new nested triangulation based on ${\bf r} \in \overline{\bf Q_{1}}$ is the one that these markings create.
In fact this is indeed the case. 
Here is the argument. As a matter of nomenclature, we will use the same concept of {\bf level} for the new triangulation $\cT^{+}({\bf r})$  
as we did before. The smallest triangles are said to be of level $0$ with side lengths equal to $\sqrt{3}$, and subsequent
sizes have levels $1,2, \dots$ of side lengths $2\sqrt{3},4\sqrt{3},\dots$.

Take any edge $e$ of a level $k+1$ triangle $T^{+}$ from the triangulation
$\cT^{+}({\bf r})$ and 
let $z$ be its midpoint. The two ends of $e$ have the form $z \pm 2^{k}3w$ for $w$ in
one of the w-directions of the lattice, see the 
red line segment in Fig.~\ref{rightBisectors}. Since both end points lie in $r_{1} + \cdots + r_{k+1}+ 2^{k+1}3 P$, we see that
$z\in r_{1}+r_{2} + \cdots + r_{k}+2^{k}3 P \subset r_{1}+ \cdots + r_{k}+2^{k}Q
=  q_{1}+ \cdots q_{k}+2^{k}Q$, and $r_{k+1}\equiv 2^{k}3w \mod 2^{k+1}3P$. Let $a$ be in the a-direction
at right-angles to $w$ and consider the two points $z \pm 2^{k}a$. These are two points
of a level $k+1$ triangle of $\cT({\bf q})$ and $z$ is its midpoint. 

To see this explicitly we take the case where $w=w_{1}$ and $a= a_{2}=  2w_{2}-w_{1}$. Then 

\[2^{k}3w - 2^{k}a = 2^{k}(3w_{1} -2w_{2} + w_{1} ) = 2^{k+1}(2w_{1}-w_{2}) \equiv 0 \mod 2^{k+1}Q\,.\]

This shows that 
\begin{align*}z+2^{k}a &\in r_{1} + \cdots + r_{k}+2^{k}3w+ 2^{k+1}Q\\
&= r_{1} + \cdots + r_{k+1}+ 2^{k+1}Q = q_{1} + \cdots + q_{k+1}+ 2^{k+1}Q\,.
\end{align*}
See Fig.~\ref{rightBisectors}. This proves:

\begin{figure}
\centering
\includegraphics[scale=0.6]{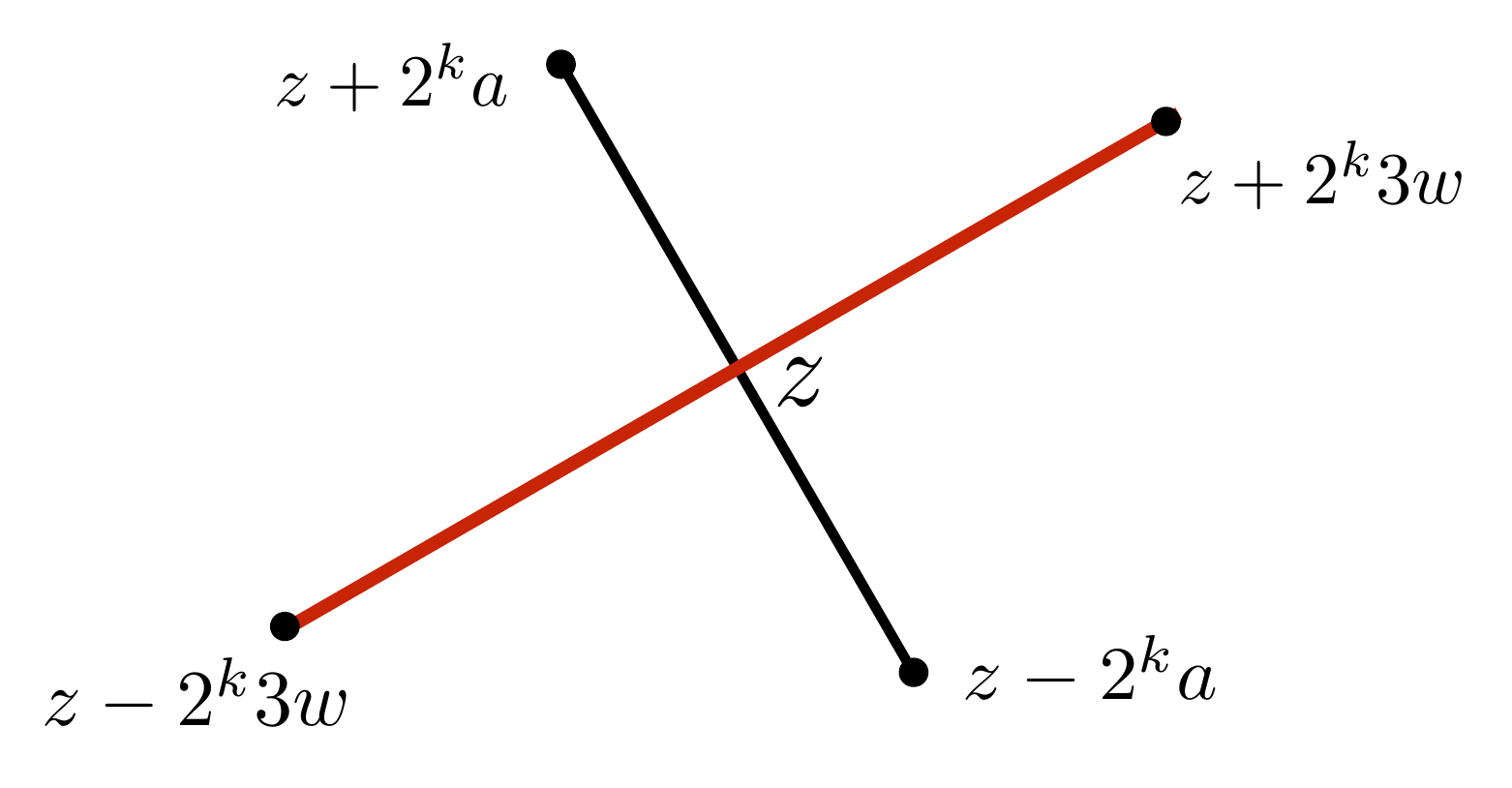}
\caption{How triangle sides of triangles from $\cT({\bf q})$ and $\cT^{+}({\bf r})$ right-bisect each other.}
\label{rightBisectors}
\end{figure}

\begin{prop}\label{rightBisectorThm} At their midpoints, the edges of level $k+1$ triangles of $\cT^{+}({\bf r})$ both
right-bisect and are right-bisected by edges of level $k+1$ of $\cT({\bf q})$.
\end{prop}

\begin{figure}
\centering
\includegraphics[scale=0.4]{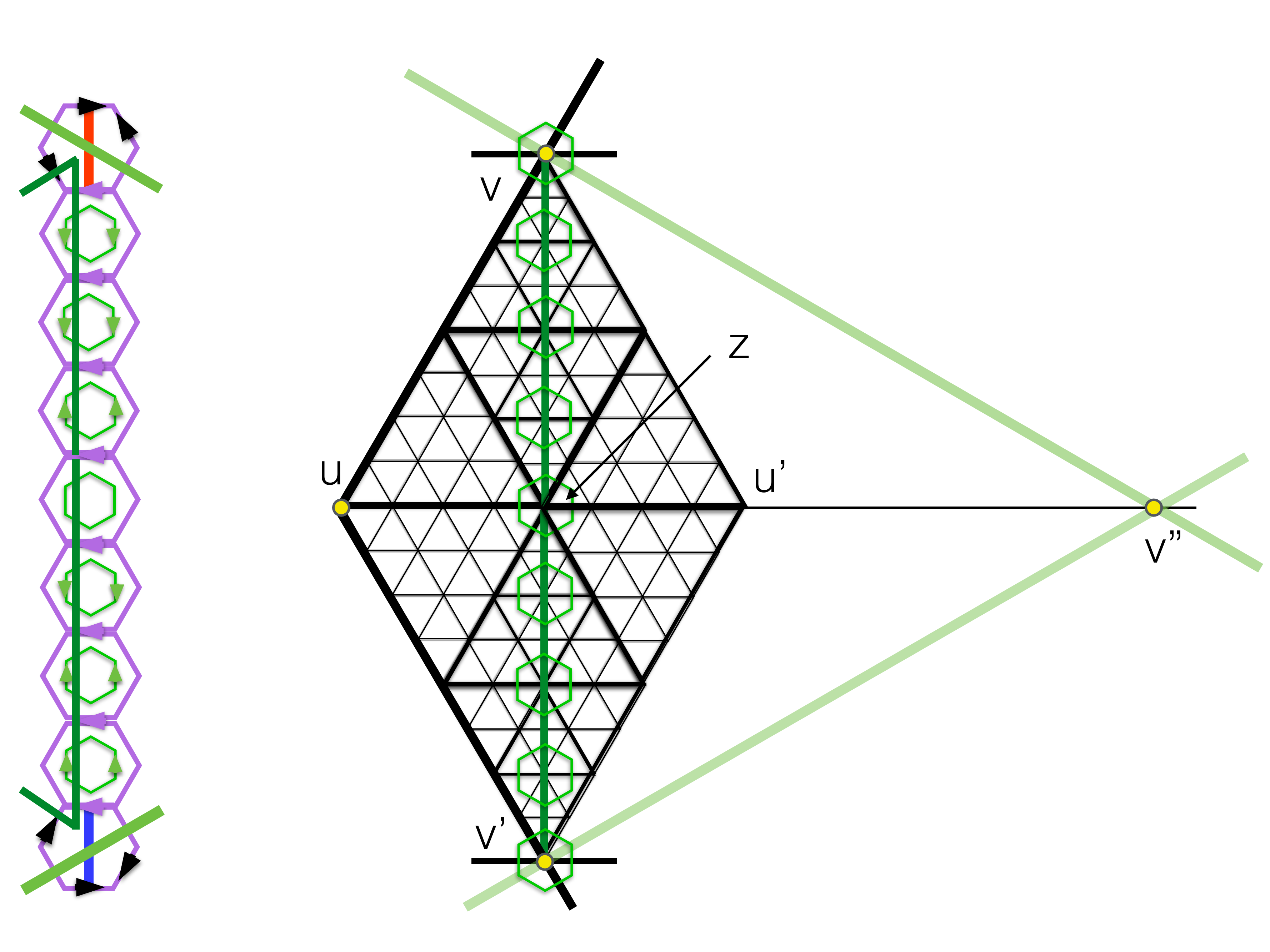}
\caption{$v,v'$ are the ends of an edge $e$ of a triangle $T^{+}$ from the nesting determined by
$\cT^{+}({\bf r})$. At its midpoint $z$ we see the edge $u, u'$ of a triangle $T$ from $\cT({\bf q})$.
The black triangles all come from the triangulation of $\cT({\bf q})$, the largest ones being of level $3$. 
The edge $e$ is maximal, in the sense that it is not part of an edge of some larger triangle from $\cT^{+}({\bf r})$.
Thus at its ends, the stripes of the large hexagons at $v$ and $v'$ are in the directions of the other sides
of $T^{+}$. The inner hexagons
along $vv'$ are centered at double hexagon centers and their stripes are 
all oriented in the same direction, namely perpendicular to the $vv'$.  At the left we have separated out the
outer hexagons that overlay the small hexagons along $vv'$.  We see their matching arrows and how their
stripes align to form
the edge $vv'$ (in green). The colors (red/blue) of the short diameters of these large hexagons
are determined by (or determine, whichever way one wants to put it) the color rule that we see in Fig.~\ref{bigHexLineUp},
though note that the stripe of
the large hexagons is perpendicular to that of the small ones, so the right/left crossing rule is opposite!
The fact that the stripe orientation changes at the end dictates that the edge $vv'$ is an interior
edge of a larger triangle. The shift indicated by the orientation of the arrows matches this.}
\label{bisectionProperty}
\end{figure}

These midpoints are centers of double hexagon tiles and, as in all double hexagon tiles, the
stripes of the inner and outer hexagons are at right-angles. This applies
to triangles of all levels $k=0,1,2, \dots $. Since by Prop.~\ref{centralEdgeProperty} every stripe of a small 
hexagon lies at the middle of some edge of some triangle, 
we conclude that the nested triangulation $\cT^{+}({\bf r})$ is directly related to the stripes on the 
large hexagon tiles, namely the path created by these stripes form the triangles of this triangulation. 
Thus the triangular tiling produced by the outer hexagons is indeed the one produced
by the nested triangulation of $\cT^{+}({\bf r})$.  Fig.~\ref{bisectionProperty} illustrates
what is going on here and also shows that the edge shifting (involved in truly nesting the 
triangles) is also properly indicated by the outer hexagon tiles.

\medskip

A direct consequence of Prop.~\ref{rightBisectorThm} is
that generic (resp. singular) $3P$-nested triangulations give rise to generic (resp. singular) $Q$-nested triangulations:

\begin{prop} \label{generic-to-generic}
Let ${\bf r} \in \overline{\bf Q_{1}}$ and let ${\bf q}$ be its image in
$\overline{\bf Q_{0}}$. Then $\cT({\bf q})$ is generic if and only if  $\cT^{+}({\bf r})$ is generic.
\qed
\end{prop}

\iffalse
\begin{prop} \label{generic-to-generic}
\textcolor{red}{Let ${\bf r} \in \overline{\bf Q_{1}}$. Then $\cT^{+}({\bf r})$ has an infinite line of $a$-direction (resp. $w$-direction) if and only if
either ${\bf r} \in Q + \overline{\Z_{2}} 3w$ for some $w\in \{w_{1}, w_{2}, w_{2}-w_{1}\}$
 {\rm (}resp.
${\bf r} \in Q + \overline{\Z_{2}} 3a$ for some $a\in\{a_{1},a_{2},a_{1}+a_{2}\}${\rm )}. 
In particular,
$\cT({\bf q})$ is generic if and only if  $\cT^{+}({\bf r})$ is generic, where the pair $({\bf q}, c)$ determines ${\bf r}$.}
\qed
\end{prop}

\fi

%%%%%%%%%%%%%%%%%%%%%%%%%%%%%%
\section{From nested triangulations to  double hexagon tilings}\label{nt2dht}
%%%%%%%%%%%%%%%%%%%%%%%%%%%%%%%

\begin{figure}
\centering
\includegraphics[scale=0.25]{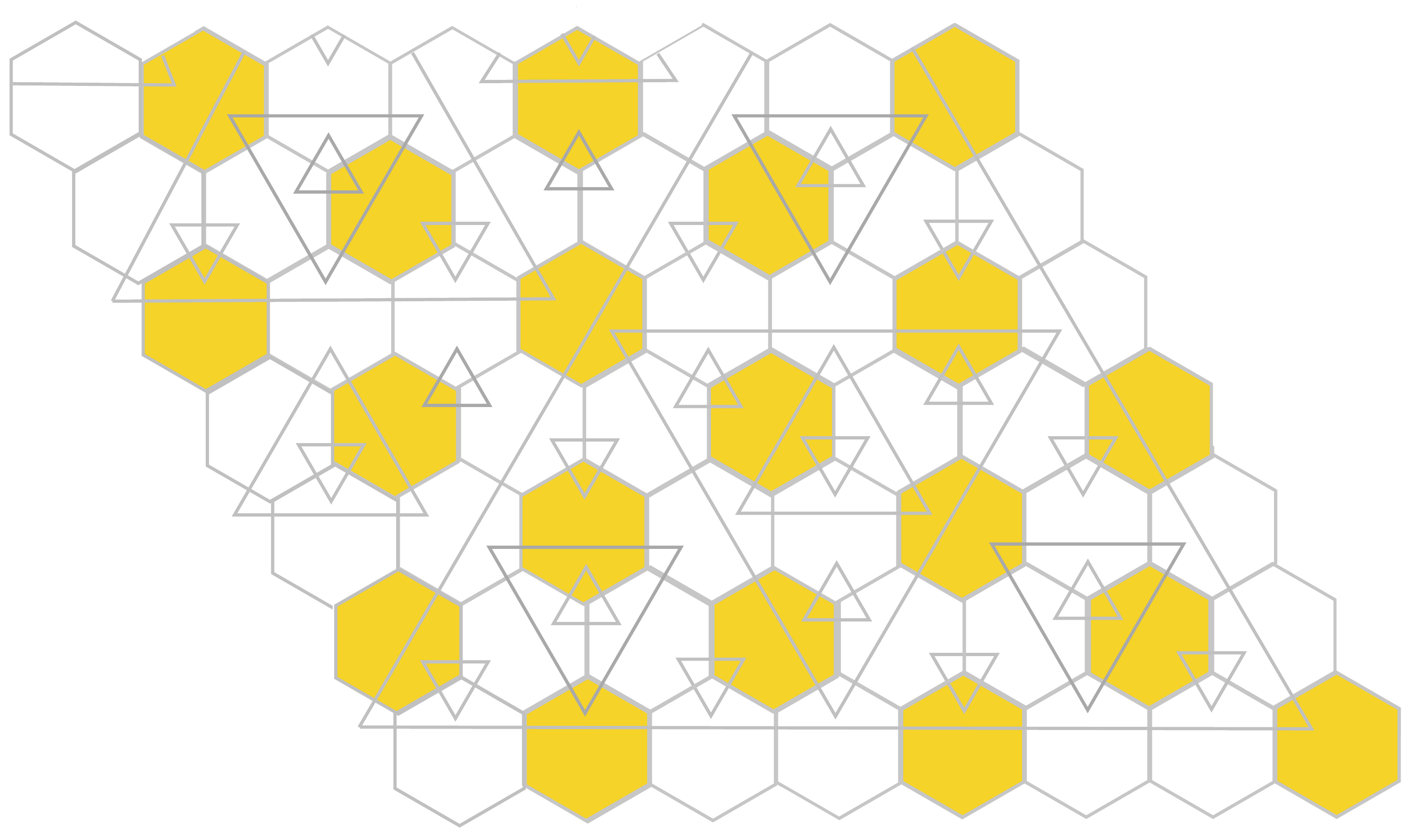}
\caption{Starting from the pair $({\bf q},{\bf r})$, with $r_{1} \equiv c \mod 3P$, hexagons centered at the lattice points of $Q$ are shown, with those centered on a coset $c +3P$ 
indicated in yellow. The nesting of the triangulation $\cT({\bf q})$ arising from ${\bf q}$
 is indicated.} 
\label{nestedTrianglesAndHexagons}
\end{figure}

At this point it is rather clear that given any triangulation $\cT({\bf q})$ and any choice of
one coset $c+ 3P$ leading to $\cT^{+}({\bf r})$ ought to lead to a legal double hexagon tiling.
Here are the details. We will assume that both triangulations are non-singular. This guarantees
that there are no infinite edges, we get proper triangulations, and they are nested. This nesting
can be geometrically manifested by laterally shifting the edges as indicated by the nesting.
The Voronoi cells of the lattices (nearest neighbor cells) are hexagons centered respectively on the points
of $\cT({\bf q})$ and $\cT^{+}({\bf r})$. We  know that every hexagon from $\cT({\bf q})$ has a triangle edge
passing through it and this edge will be shifted laterally in nesting. This is shown in
Fig.~\ref{nestedTrianglesAndHexagons}. The hexagon is made into a well-arrowed
hexagon by placing the pair of parallel arrows in the direction of the shift. The small hexagons now make
a well-arrowed and well-matched hexagon tiling. 

\begin{figure}[h]
\centering
\includegraphics[scale=0.25]{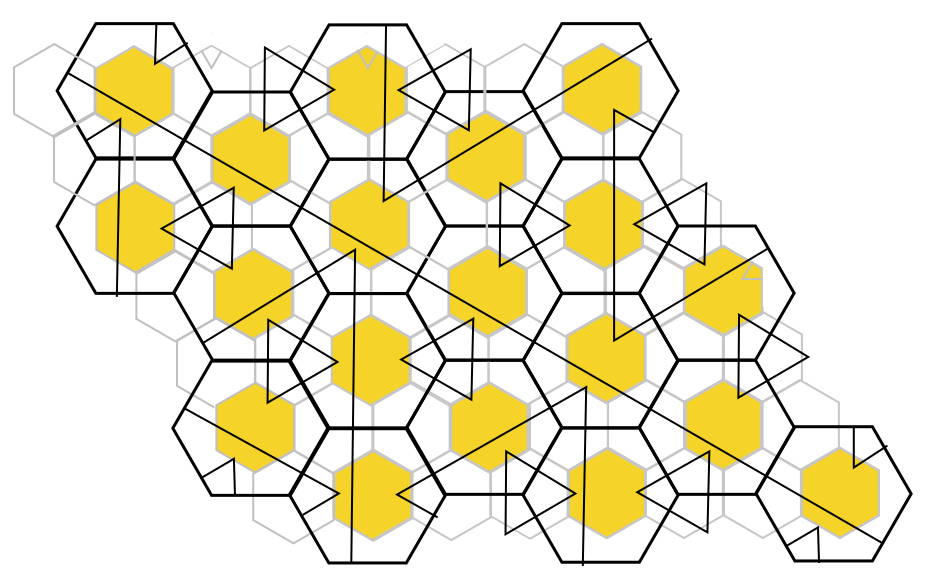}
\caption{Continuing from Fig.~\ref{nestedTrianglesAndHexagons}, from ${\bf r} \in \overline{\bf Q_{1}} $ we obtain a nested triangulation $\cT^{+}({\bf r})$, part of which is shown here. }
\label{nestedTrianglesAndHexagons-II}
\end{figure}

Now we do the same thing with the triangulation $\cT^{+}({\bf r})$, leading to the new hexagonal tiling, 
again with arrows indicating edge shifting in the nesting, see Fig.~\ref{nestedTrianglesAndHexagons-II}. This is the second well-arrowed and well-matched
hexagonal tiling. The outcome is that we have a tiling of double hexagon tiles which is non-singular
and legal, see Fig.~\ref{doubleHexPatch}.

%%%%%%%%%%%%%%%%%%%%%%%%%%%%%%%%%%%%%%%%%%
\section{Penrose tilings, Taylor-Socolar tilings, and beyond}
%%%%%%%%%%%%%%%%%%%%%%%%%%%%%%%%%%%%%%%%%%

\begin{figure}\centering
\includegraphics[scale=0.3]{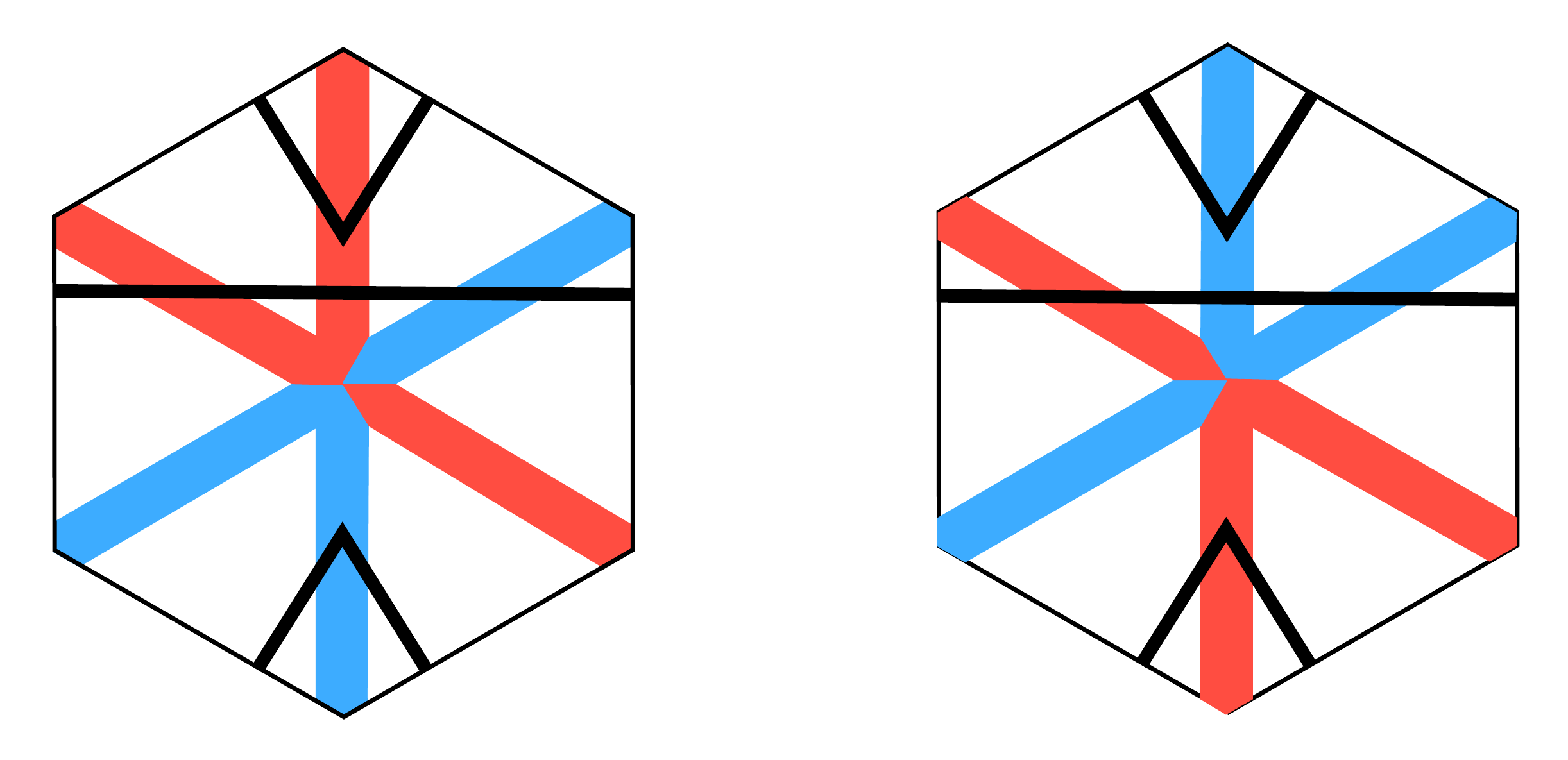}
\caption{The two prototiles for the Taylor-Socolar tilings}
\label{basic-T-Stiles}
\end{figure}

By definition a Penrose tiling is precisely a legal double hexagon tiling. Taylor-Socolar tilings (T-S tilings)
are usually defined by 
the T-S tiles shown in Fig. \ref{basic-T-Stiles}
and they are assembled as regular hexagonal tilings, but under the matching rules 
\begin{itemize}
\item[{\bf RT1}] the black lines must join continuously when tiles abut;
\item[{\bf RT2}] the ends of the diameters of two hexagonal tiles that are separated by an edge of another tile must be of opposite colors.
\end{itemize}

In Fig.~\ref{bigHexLineUp} we have seen that the diagonals of the inner hexagons of a double hexagon tiling can
be colored, and if we restrict this coloring to the actual physical area of the inner hexagons then 
we have the colorings of Fig.~\ref{basic-T-Stiles}. The matching of the outer arrows of the 
double hexagon tiling amounts to the color rule {\bf RT2}, so we have in this way one third
of a T-S tiling. If the double hexagon tiling is legal then we know that this partial hexagonal tiling of
inner hexagons along with the corner hexagons completes to a new properly arrowed hexagonal tiling
together with a corresponding nested triangulation. If we assume that the nested triangulation is non-singular, which
is generically the case, then this tiling-triangulation corresponds to a unique T-S
tiling. That is, the one-third tiling we have completes uniquely to a T-S tiling. The proof of this is
given in \cite{TSMS}---each non-singular nested triangulation corresponds to a unique T-S tiling and vice-versa. 

Thus every non-singular Penrose tiling produces inside it a non-singular T-S tiling made out of 
its inner and corner hexagons. 
Now let us go in the other direction. If we begin with a non-singular T-S tiling then it produces a nested
triangulation out of the stripes on each hexagon. Relative to a fixed coordinate system, this triangulation 
corresponds to an element ${\bf q} $ in the  $Q$-adic completion $\overline{\bf Q_{0}}$. In order to obtain a double hexagon tiling from this
we need to select which hexagons will be the inner hexagons and which the corner hexagons for the new tiling.
This amounts to choosing one coset from the $Q/3P$.
Choose one, say,  $c+ 3P$. Then there is a unique ${\bf r} \in \overline{\bf Q_{1}}$ that maps ${\bf q}$ under the natural mapping
of $\overline{\bf Q_{1}}$ to $\overline{\bf Q_{0}}$ and for which $r_{1} \equiv c \mod 3P$. This produces the centers
and the nested triangulation that determines a
legal double hexagon tiling, as we have pointed out in \S\ref{nt2dht}. 

To reiterate, we see that the nested triangulation of the large hexagon tiles determines
the nested triangulation of the inner hexagonal tiles. Thus, although we only see the coloring of one third
of an underlying T-S tiling, the entire nesting of the triangulation arising from the small hexagon tiles 
is implicitly known from the
nested triangulation of the larger hexagon tiles: we know ${\bf q}$ once we know ${\bf r}$.

A noteworthy observation comes by comparing Fig.~\ref{basic-T-Stiles}
and Fig.~\ref{NewDeco-on-PenroseTile}: it shows that the distinction between the parity
(that is, the difference between the two types of small hexagon tiles (respectively large hexagon tiles) is the
same for the T-S tiling and the Penrose tiling. Thus the parity distribution of a Penrose tiling
is the same as the parity distribution of one coset modulo $3P$ of the T-S tiles.

Although it is shown in \cite{BGG} that the two tiling spaces generated by T-S tilings and Penrose tilings define distinct MLD (mutual local derivabillity) classes, it is clear by now that the two types of tilings are intimately related,
and indeed, modulo the choice of a coset, there is a mutual derivability. 
We can summarize some key points as follows: 
\begin{theorem}
\begin{itemize}
\item[{\rm (i)}] Taylor-Socolar tilings and Penrose tilings are aperiodic.
\item[{\rm (ii)}] Given a non-singular Taylor-Socolar tiling on $Q$, one can build, in a canonical way, three different non-singular Penrose tilings, one for each of the three cosets $c+ 3P$ of $Q\mod 3P$.  At any point of $c+ 3P$ one knows exactly what type of Penrose tile should be put in that position, and this uses only local information of the T-S tiling.
\item[{\rm (iii)}] Given a non-singular Penrose tiling on some coset $c+3P$, there is a unique nested triangulation on $Q$ formed by the decoration of inner hexagons and corner hexagons of Penrose tiles. This nested triangulation gives a unique non-singular T-S tiling. Note that unlike the situation in {\rm (ii)}, this construction, is not local.
\end{itemize}
\end{theorem}

The process of producing double hexagon tiles from a pair $({\bf q}$, ${\bf r})$ suggests that we might do it again, choosing a coset $d + 3Q$ with $d \equiv c \mod 3P$ and then determining ${\bf s} \in \overline{\bf Q_{2}}$.
 This triple $({\bf q}, {\bf r},{\bf s})$ leads to triple-hexagon tiles and a triple-hexagon tiling. 
 The rules for admissibility follow the same principles as we have used above. The largest hexagonal tiles
 have middle sized hexagonal tiles at their centers, and create middle sized corner tiles around them. The
 requirement is well-arrowing throughout. This yields a well-arrowed hexagonal tiling of middle sized tiles. In the 
 same way, create from this a hexagonal tiling of small tiles, and where again we require well-arrowing throughout. 
 
\begin{figure}
\centering
\includegraphics[scale=0.4]{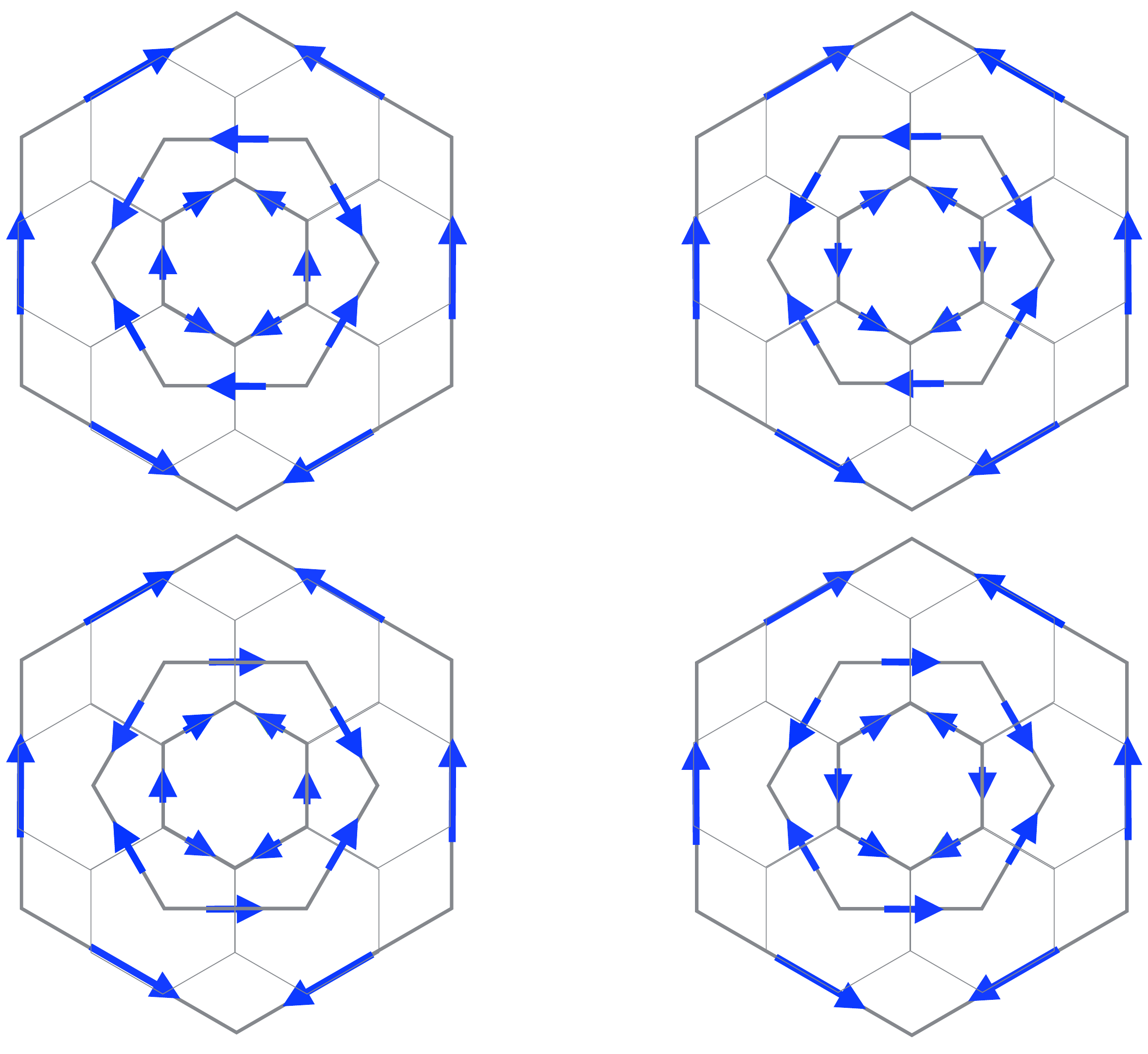}
\caption{The four types of triple (or $3$-tiered) hexagons}
\label{3-tiered-hexagons}
\end{figure}

 These triple tiles come in four types and produce a new type of hexagonal tiling, Fig. \ref{3-tiered-hexagons}. 
 There is no reason to stop there. This new hierarchical situation is illustrated in the commutative diagram \eqref{commutativeDiagram2}.
 
\begin{align}\label{commutativeDiagram2}
\overline{\bf Q_{0}}: &\quad &Q/Q  & \quad \longleftarrow & Q/2Q & \quad \longleftarrow & Q/4Q & \quad \longleftarrow &   Q/8Q & \leftarrow\cdots  \nonumber\\
\uparrow &&\uparrow & & \uparrow & & \uparrow & &  \uparrow & &\nonumber\\
\overline {\bf Q_{1}}: & \quad &Q/3P  & \quad \longleftarrow & Q/6P & \quad \longleftarrow& Q/12P & \quad
\longleftarrow & Q/24P &\longleftarrow \cdots \nonumber \\
\uparrow &&\uparrow & & \uparrow & & \uparrow & &  \uparrow & &\\
\overline {\bf Q_{2}}: & \quad &Q/3Q  & \quad \longleftarrow & Q/6Q & \quad \longleftarrow& Q/12Q & \quad
\longleftarrow & Q/24Q &\longleftarrow \cdots \nonumber \\
\uparrow &&\uparrow & & \uparrow & & \uparrow & &  \uparrow & &\nonumber\\
\overline {\bf Q_{3}}: & \quad &Q/9P  & \quad \longleftarrow & Q/18P & \quad \longleftarrow& Q/36P & \quad
\longleftarrow & Q/72P &\longleftarrow \cdots \nonumber \\
\uparrow &&\uparrow & & \uparrow & & \uparrow & &  \uparrow & &\nonumber\\
\vdots &&\vdots & & \vdots & & \vdots  & &  \vdots  & &\nonumber\\
\nonumber
\end{align}

More generally there are `$n$-tuple hexagons' or to have a 
better sounding name, $n$-{\bf tiered hexagons}, each of which consists of $n$
hexagons stacked within each other, which tile the plane according to 
the $n$-th line of \eqref{commutativeDiagram2}. There are two choices
for the orientation of a hexagon at each stage of layering, but after taking into
account rotations, this gives $2^{n-1}$ types of these tiered tiles.
Non-singularity (resp. singularity) is a common property to all levels.

\section{Outlook}\label{outlook}
The purpose of this paper has been to clarify the unity
that exists between the Taylor-Socolar tilings and the Penrose
hexagonal tilings---a unity that can be expressed both geometrically
and algebraically in terms of double hexagon tiles. Each non-singular
legal double hexagon tiling encompasses both a Penrose tiling and
a T-S tiling, and this pairing can be interpreted algebraically in terms
of \eqref{commutativeDiagramSmall}.  Each of the two hexagonal tilings
leads to a nested triangulation, and these two are bound together by the 
simple rule that triangle edges of each right-bisect edges of the other. 

There are two issues that arise here that we have not discussed, but 
plan to pursue in a future work. The first is the nature of singularities
in these tilings from both the geometric and algebraic perspectives, and
their detailed manifestation in the corresponding tiling hulls. 
The second is the study of the $n$-tiered hexagonal tilings. The algebraic
setting which uses the first $n$ rows of the commutative diagram
\eqref{commutativeDiagram2} suggests that the $n$-tiered hexagons
lead to aperiodic tilings in which there are potentially $2^{n-1}$ types
of tiles. Thus there is a hierarchy of aperiodic hexagonal tilings, and their
corresponding tiling hulls, about which we know very little. 
%%%%%%%%%%%%%%%%%%%%%

\section*{Acknowledgment}

The first author would like to acknowledge
that this work was supported by a National Research Foundation
of Korea (NRF) Grant funded by the Korean
Government (MSIP) (No. 2014004168) and by research fund of Catholic Kwandong University(CKURF-201604560001). She is also grateful for the support by the Korea
Institute for Advanced Study (KIAS) .
%%%%%%%%%%%%%%%%%%%%%

\end{document}